\documentclass[11pt]{article}
\usepackage{graphicx}
\usepackage{amsmath}
\usepackage{amssymb}
\usepackage{amsthm}
\usepackage{url}
\newcommand\nonu{\nonumber}
\newcommand\sa{\smallskipamount}
\newcommand\ba{\bigskipamount}
\newcommand\bLP{\\[\ba]}
\newcommand\sPP{\\[\sa]\indent}
\newcommand\bPP{\\[\ba]\indent}
\newcommand\RR{\mathbb{R}}
\newcommand\ZZ{\mathbb{Z}}
\newcommand\FSJ{{\cal J}}
\newcommand\FSK{{\cal K}}
\newcommand\FSO{{\cal O}}
\newcommand\FSP{{\cal P}}
\newcommand\al\alpha
\newcommand\be\beta
\newcommand\de\delta
\newcommand\ep\varepsilon
\newcommand\tha\theta
\newcommand\ka\kappa
\newcommand\la\lambda
\newcommand\om\omega
\newcommand\Ga{\Gamma}
\newcommand\De{\Delta}
\newcommand\half{\frac12}
\newcommand\thalf{\tfrac12}
\newcommand\iy\infty
\newcommand\pa\partial
\newcommand\Zpos{\ZZ_{>0}}
\newcommand\Znonneg{\ZZ_{\ge0}}
\newcommand\const{{\rm const.}\,}
\newcommand{\hyp}[5]{\,\mbox{}_{#1}F_{#2}\!\left(
  \genfrac{}{}{0pt}{}{#3}{#4};#5\right)}
\newcommand\LHS{left-hand side}
\newcommand\RHS{right-hand side}
\newcommand\supp{\mathop{\rm supp}\nolimits}

\numberwithin{equation}{section}
\newtheorem{theorem}{Theorem}[section]
\newtheorem{proposition}[theorem]{Proposition}
\newtheorem{lemma}[theorem]{Lemma}
\newtheorem{Remark}[theorem]{Remark}
\newenvironment{remark}{\begin{Remark}\rm}{\end{Remark}}
\newcommand\Proof{\noindent{\bf Proof}\quad}

\begin{document}
\title{Differentiation by integration using orthogonal polynomials,\\
a survey}
\author{Enno Diekema\footnote{\tt e.diekema@gmail.com}\;
and Tom H. Koornwinder\footnote{Korteweg-de Vries Institute,
University of Amsterdam, \tt thkmath@xs4all.nl}}
\date{}
\maketitle
\begin{abstract}
This survey paper discusses the history of approximation formulas for
$n$-th order derivatives by integrals involving orthogonal polynomials.
There is a large but rather disconnected corpus of literature on such
formulas. We give some results in greater generality than in the literature.
Notably we unify the continuous and discrete case.
We make many side remarks, for instance on wavelets,
Mantica's Fourier-Bessel functions and Greville's minimum $R_\al$ formulas
in connection with discrete smoothing.
\end{abstract}

%
%
\section{Introduction}
In many applications one needs to estimate or approximate the first or
higher derivative of a function which is only given in sampled form
or which is perturbed by noise. Good candidates for an approximation
of the first derivative $f'(x)$ are the two expressions
\begin{align}
&\frac3{2\delta}\,\int_{-1}^1 f(x+\de\xi)\,\xi\,d\xi
\label{110}\\
\noalign{\mbox{and}}
&\frac3{2N(N+\thalf)(N+1)\de}\sum_{\xi=-N}^N f(x+\de\xi)\,\xi
\label{111}
\end{align}
for $\de$ small. The first one is continuous, the second one discrete.
These formulas have a long history going back to
Cioranescu \cite{31} (1938),
Haslam-Jones \cite{26} (1953),
Lanczos \cite[(5-9.1)]{19} (1956) and
Savitzky \& Golay \cite{10} (1964).
What remains hidden in \eqref{110} and \eqref{111} is that the factor $\xi$
in the integrand or summand can better be considered as an
orthogonal polynomial of degree 1 with respect to a constant weight function
on $[-1,1]$ (in case of \eqref{110}) or with respect to constant weights
on $\{-N,-N+1,\ldots,N\}$ (in case of \eqref{111}). With this point of view
it can be immediately shown that \eqref{110} and \eqref{111} tend to $f'(x)$
as $\de\downarrow0$. Moreover the way is opened to a far reaching
generalization of \eqref{110} and \eqref{111} for the approximation
of higher derivatives and with the involvement of
general orthogonal polynomials.
Such results were already given by Cioranescu \cite{31} in 1938.

Curiously enough, none of the later papers mentioned above is referring
to one of the earlier papers. The results of
Cioranescu \cite{31} and Haslam-Jones \cite{26}
were hardly taken up by anybody. On the other hand
Lanczos \cite[(5-9.1)]{19} and in particular Savitzky \& Golay \cite{10}
had a lot of follow-up by others. One reason for this citation behaviour is
probably that Cioranescu and Haslam-Jones were pure analysts,
Lanczos was an applied mathematician working in numerical analysis,
and Savitzky \& Golay were motivated by spectroscopy, considered as a part
of chemistry.

It is the aim of the present paper to give a survey of these results,
developments and further considerations suggested by them. Moreover,
we formulate some results in a more general way than has probably appeared
before in literature. It was for us a surprise to see that so many
different parts of classical analysis and of applied mathematics are tied
together by this theme. All papers until now only treated smaller parts of
this wide field. We hope to share with our readers the pleasure to have
a comprehensive view.

The present work stems from a long practice by the first author
in signal processing in applied situations, where he met the problem 
of differentiating an analog signal (and later a sampled signal) disturbed by
noise, see for instance Strik \cite{9}.
The main problem occurring there was the difficulty
to build (or program) an ideal differentiator,
because the noise of the system will
cause an instability. Therefore an integrating factor for the high
frequencies to the differentiator is needed.
When the signal is sampled the same problem occurs (Hamming \cite{7}).
Without being aware of the literature mentioned in the beginning
of this Introduction, he tried to use an integrating factor
by the method of the least squares and then he independently found 
special cases of the approximation formulas for higher derivatives
by integrals involving all classical orthogonal polynomials
as well as the Chebyshev polynomials of a discrete variable.
He never published the results, but he used them as material for a course
in stochastic system theory at the "Saxion Hogeschool" in Enschede.

The contents of this paper are as follows. In section 2 we give preliminaries
on orthogonal polynomials and on the Taylor formula.
In section 3 we start formulating the approximation theorem in great generality
and next discuss how the contributions of Cioranescu, Haslam-Jones and Lanczos
are related to this general theorem. We emphasize the important role of
least-square approximation behind this theory. 
Our discussion gives room for several side observations, for instance
on Jacobi type polynomials and on wavelets.
Section 4 is focused on the discrete case and the applications to filters.
We start with a multi-term extension of the main theorem in section 3.
Its special case of constant weights contains the seminal work of
Savitzky \& Golay. We introduce the characteristic (or transfer) function
and we make connection with Mantica's \cite{62} Fourier-Bessel functions.
In the smoothing case we discuss the work by Greville \cite{45} (based on
older work by Sheppard \cite{63}) on so-called minimum $R_\al$ and
minimum $R_\iy$ formulas.
In the Appendix we give new derivations of the
characteristic functions for these cases,
using Hahn and Krawtchouk polynomials.
The $R_\iy$ case connects with another survey paper \cite{47} by the
second author and Schlosser.
Finally, in section 5, we discuss log-log plots of transfer functions in
some simple cases.
\section{Preliminaries}
\subsection{Orthogonal polynomials}
\label{18}
Let $\mu$ be a positive Borel measure on $\RR$ with infinite support
(or equivalently a nondecreasing function on $\RR$ with an infinite number
of points of increase) such that
$\int_\RR |x|^n\,d\mu(x)<\iy$ for all $n\in\Znonneg$.
Consider polynomials $p_n$ ($n\in\Znonneg$) of degree $n$ such that
\begin{equation}
\int_\RR p_m(x)\,p_n(x)\,d\mu(x)=0\qquad(m\ne n).
\label{1}
\end{equation}
The $p_n$ are {\em orthogonal polynomials} with respect to
the measure $\mu$, see for instance Szeg{\H{o}} \cite{1}. Up to
constant nonzero factors they are uniquely determined by the above
properties.
If $\mu$ has support within some closed interval $I$ then we say that
the $p_n$ are orthogonal polynomials with respect to $\mu$ on $I$.
Typical cases of the orthogonality measure $\mu$ are:
\begin{enumerate}
\item
$d\mu(x)=w(x)\,dx$ on $I$ with the {\em weight function} $w$ a
nonnegative integrable function on $I$. Then \eqref{1} takes the form
\[
\int_I p_m(x)\,p_n(x)\,w(x)\,dx=0\qquad(m\ne n).
\]
\item
$\mu$ has discrete infinite support $\{x_0,x_1,x_2,\ldots\}$. So there
are positive numbers $w_0,w_1,w_2,\ldots$ ({\em weights}) such that
\eqref{1} takes the form
\begin{equation}
\sum_{k=0}^\iy p_m(x_k)\,p_n(x_k)\,w_k=0\qquad(m\ne n).
\label{87}
\end{equation}
\item
Contrary to what was supposed earlier, we can also consider the case that
$\mu$ has finite support $\{x_0,x_1,\ldots,x_N\}$ with corresponding
weights $w_0,w_1,\ldots,w_N$. Then we have orthogonal polynomials
$p_n$ only for $n=0,1,\ldots,N$ and \eqref{1} takes the form
\begin{equation}
\sum_{k=0}^N p_m(x_k)\,p_n(x_k)\,w_k=0\qquad(m\ne n).
\label{88}
\end{equation}
\end{enumerate}

Special examples of case 1 are given by the
{\em classical orthogonal polynomials}
(Jacobi, Laguerre and Hermite polynomials).
In particular, we will meet
the {\em Legendre polynomials} $P_n$, which are special Jacobi polynomials
and where $I=[-1,1]$, $w(x)=1$ and $P_n(1)=1$.

A special example of case 3 are the  {\em Hahn polynomials}
$x\mapsto Q_n(x;\al,\be,N)$
for $\al=\be=0$ ($n=0,1,\ldots,N$). Here $x_i=i, w_i=1$
($i=0,1,\ldots,N$) and $Q_n(0;0,0,N)=1$ (see \cite[\S9.5]{5} and
references given there, or \cite[Ch.~2]{3}, where another notation is
used). Hahn polynomials of general parameters were already introduced
in 1875 by Chebyshev \cite{18},
long before Hahn, but the above special case of constant weights is in
particular named after Chebyshev, although
in a slightly different notation and normalization.
See {\em Chebyshev's polynomials of a discrete variable} $t_n(x)=t_n(x,N)$
($n=0,1,\ldots,N-1$) in \cite[\S2.8]{1}, \cite[\S10.23]{4}.
They are orthogonal on the set $\{0,1,\ldots,N-1\}$
with respect to constant weights 1. Hence we must have that
$t_n(x,N)=\const\,Q_n(x;0,0,N-1)$.
The constant can be computed by comparing the recurrence relation
\cite[(9.5.3)]{5} for $\al=\be=0$ and $N$ replaced by $N-1$ with the
recurrence relation \cite[10.23(6)]{4}. Then we obtain:
\begin{equation}
t_n(x,N)=(-N+1)_n\,Q_n(x;0,0,N-1),
\label{11}
\end{equation}
where $(a)_n:=a(a+1)\ldots(a+n-1)$ is the
{\em Pochhammer symbol}. Thus $t_n(0,N)=(-N+1)_n$.
These polynomials are also known as {\em Gram polynomials},
see \cite[\S7.13 and \S7.16]{16}. This last name we will use in this paper.
The polynomials $x\mapsto t_n(x,N)$ have the shifted Legendre polynomials
$x\mapsto P_n(2x-1)$ (orthogonal on $[0,1]$ with respect to a constant
weight function) as a limit case (see \cite[(2.8.6)]{1}):
\begin{equation}
\lim_{N\to \infty}N^{-n}\,t_n(Nx,N)=P_n(2x-1)
\label{45}
\end{equation}

For given orthogonal polynomials $p_n$
define the constants $h_n$ and $k_n$ by
\begin{equation}
h_n:=\int_\RR p_n(x)^2\,d\mu(x),\qquad
p_n(x)=k_nx^n+\mbox{terms of degree less than $n$}.
\label{10}
\end{equation}
\begin{lemma}
We have
\begin{equation}
\int_\RR p_n(x)\,x^n\,d\mu(x)=\frac{h_n}{k_n}\
\label{2}
\end{equation}
\end{lemma}
\Proof
We have $k_nx^n=p_n(x)+q_{n-1}(x)$ with $q_{n-1}$ a polynomial of
degree $<n$. Hence
\[
k_n\int_\RR p_n(x)\,x^n\,d\mu(x)=
\int_\RR p_n(x)^2\,d\mu(x)+\int_\RR p_n(x)\,q_{n-1}(x)\,d\mu(x)=
h_n+0=h_n
\eqno{\qed}
\]

One of the properties which characterize the classical
orthogonal polynomials is that they are given by a (generalized)
{\em Rodrigues formula}
\begin{equation}
p_n(x)=\frac1{K_n\,w(x)}\,\frac{d^n}{dx^n}\,(w(x)\,X(x)^n)
\label{103}
\end{equation}
(see \cite[10.6(1)]{4}).
Here $X$ is a polynomial of degree $\le2$ and
\begin{equation}
K_n=\frac{(-1)^n\,k_n\,n!}{h_n}\,\int_\RR(X(x))^n\,d\mu(x).
\label{104}
\end{equation}
For the proof of \eqref{104} substitute \eqref{103} and $d\mu(x)=w(x)\,dx$
in \eqref{2} and perform integration by parts $n$ times.

The explicit values of $h_n$ and $k_n$ defined by \eqref{10} can be
given in our two main examples:
\begin{itemize}
\item
Legendre polynomials $P_n$ (see \cite[(9.8.63), (9.8.65)]{5}):
\begin{equation}
h_n=\frac2{2n+1},\qquad
k_n=2^{-n}\binom{2n}n.
\label{13}
\end{equation}
From this we immediately obtain the values of $h_n$ and $k_n$ in
the case of shifted Legendre polynomials $x\mapsto P_n(2x-1)$:
\begin{equation}
h_n=\frac1{2n+1},\qquad
k_n=\binom{2n}n.
\label{59}
\end{equation}
\item
Gram polynomials $x\mapsto t_n(x,N)$
(combine \eqref{11} with \cite[(9.5.2), (9.5.4)]{5}):
\begin{equation}
h_n=\frac{(N-n)_{2n+1}}{2n+1}\,,\qquad
k_n=\binom{2n}n.
\label{12}
\end{equation}
From this we immediately obtain the values of $h_n$ and $k_n$ in
the case of centered Gram polynomials $x\mapsto t_n(x+N,2N+1)$:
\begin{equation}
h_n=\frac{(2N+1-n)_{2n+1}}{2n+1}\,,\qquad
k_n=\binom{2n}n.
\label{68}
\end{equation}
\end{itemize}

The {\em reproducing kernel} for the space $\FSP_n$ of polynomials
of degree $\le n$
in the Hilbert space $L^2(\RR,\mu)$ is given by
\begin{equation}
{\bf K}_n(x,y):=\sum_{k=0}^n \frac{p_k(x)\,p_k(y)}{h_k}\qquad(x,y\in\RR).
\label{31}
\end{equation}
Then (see \cite[Remark 5.2.2]{34}) the
{\em Christoffel-Darboux formula} gives
\begin{equation}
{\bf K}_n(x,y)=\frac{k_n}{k_{n+1}\,h_n}\,
\frac{p_{n+1}(x)p_n(y)-p_n(x)p_{n+1}(y)}{x-y}\quad(x\ne y).
\label{62}
\end{equation}
The integral operator $\FSK_n$ corresponding to \eqref{31} is given by
\begin{equation}
(\FSK_n f)(x):=\int_{\RR}f(y)\,K_n(x,y)\,d\mu(y)\qquad
(f\in L^2(\RR),\;x\in\RR),
\end{equation}
It is the orthogonal projection of the Hilbert space $L^2(\RR,\mu)$
onto $\FSP_n$.
In particular,
\begin{equation}
\FSK_n f=f\qquad(f\in\FSP_n).
\label{30}
\end{equation}
Furthermore, for $f\in L^2(\RR,\mu)$, $\FSK_n f$ is the element of $\FSP_n$
which is on minimal distance to $f$ (in the metric of
the Hilbert space $L^2(\RR,\mu)$).

Later we will need the following.
For Legendre polynomials formula \eqref{62} for $y=1$ and with $n$ replaced by
$n-1$ becomes:
\begin{align}
{\bf K}_{n-1}(x,1)&=\thalf n\,\frac{P_n(x)-P_{n-1}(x)}{x-1}
\nonu\\
&=\thalf n P_{n-1}^{(1,0)}(x)
\label{66}\\
&=\thalf(P_n'(x)+P_{n-1}'(x)).
\label{67}
\end{align}
The $P_{n-1}^{(1,0)}(x)$ in \eqref{66} is a Jacobi polynomial
(see \cite[Ch.~4]{1}). We used \eqref{13} in the first equality,
\cite[10.8(32)]{4} in \eqref{66}, and
\cite[10.10(13), 10.10(14)]{4} in \eqref{67}.

\subsection{Taylor formula}
Recall a version of Taylor's theorem formulated by Hardy
\cite[\S151]{30}:
\begin{proposition}
\label{56}
Let $x\in\RR$ and let $I$ be an interval containing $x$.
Let $f$ be a continuous function on $I$ such that its derivatives of
order $1,2,\ldots,n$ at $x$ exist. Then
\begin{equation}
f(y)=\sum_{k=0}^n\frac{f^{(k)}(x)}{k!}\,(y-x)^k+o(|y-x|^n)\quad
\mbox{as $y\to x$ on $I$.}
\label{75}
\end{equation}
\end{proposition}

In this proposition the derivatives should be interpreted as right or
left derivatives if $x$ is an endpoint of the interval $I$.
(Although this special case is not explicit in Hardy's formulation,
it is also a consequence of his proof.)$\;$

Proposition \ref{56} suggests a notion more general than $n$-th derivative.
Let $f$ be a continuous function on an interval $I$ containing $x$ and let
there be constants $c_0,c_1,\ldots,c_n$ such that
\begin{equation}
f(y)=\sum_{k=0}^n\frac{c_k}{k!}\,(y-x)^k+o(|y-x|^n)\quad
\mbox{as $y\to x$ on $I$.}
\label{76}
\end{equation}
Then we call $c_n$ the $n$-th {\em Peano derivative} of $f$ at $x$.
This definition goes back to Peano \cite{32} in 1891.
By Proposition \ref{56} the existence of $f^{(n)}(x)$ implies the existence
of the $n$-th Peano derivative, equal to $f^{(n)}(x)$.
The converse implication is true for $n=1$ but not for $n>1$,
see a counterexample in
\cite[Example 1.2]{33}.

For later use we restate Proposition \ref{56} as follows:
\begin{proposition}
\label{54}
With the assumptions of Proposition \ref{56} we have
\begin{equation}
f(x+\de)=\sum_{k=0}^n\frac{f^{(k)}(x)}{k!}\,\de^k+\de^n F_{x,n}(\de)
\label{57}
\end{equation}
with $F_{x,n}$ continuous on $I-x$ and $F_{x,n}(0)=0$.
Furthermore, $F_{x,n}$ is bounded on $I-x$ if $f$ is bounded on $I$.
Finally, if $I$ is unbounded and $f$ is of polynomial growth on $I$
then $F_{x,n}$ is of polynomial growth on $I-x$.
\end{proposition}

We can say more about the remainder term in \eqref{57} if moreover
$f\in C^n(I)$ (see for instance Apostol \cite[Theorem 7.6]{2}):
\begin{proposition}
\label{81}
Keep the assumptions of Proposition \ref{56}. Assume moreover that
$f\in C^n(I)$. Then for $F_{x,n}(\de)$ in \eqref{57} we have
\begin{equation}
F_{x,n}(\de)=\frac1{(n-1)!}\,
\int_0^1 \big(f^{(n)}(x+t\de)-f^{(n)}(x)\big)\,(1-t)^{n-1}\,dt
\label{82}
\end{equation}
and the function $(x,y)\mapsto F_{x,n}(y-x)$
is continuous on $I\times I$.
If $f^{(n)}$ is of polynomial growth on $I$ then $F_{x,n}(\de)\to0$
as $\de\to0$ uniformly for $x$ in compact subsets of $I$.
\end{proposition}
\section{Higher derivatives approximated by integrals}
Let us first state and prove the main theorem and next discuss the many
instances of it in the literature, usually more restricted but
occasionally more general than our formulation.
\begin{theorem}
\label{16}
For some $n$ let $p_n$ be an orthogonal polynomial of degree $n$ with
respect to the orthogonality measure $\mu$.
Let $x\in\RR$. Let $I$ be a closed interval such that, for some
$\ep>0$, $x+\delta\xi\in I$ if $0\le \delta\le\ep$ and
$\xi\in\supp(\mu)$.
Let $f$ be a continuous function on $I$ such that its derivatives of order
$1,2,\ldots,n$ at $x$ exist. In addition, if
$I$ is unbounded, assume that $f$ is of at most polynomial
growth on $I$. Then
\begin{equation}
f^{(n)}(x)=\frac{k_n n!}{h_n}\,\lim_{\delta\downarrow0}\frac1{\delta^n}\,
\int_\RR f(x+\de\xi)\,p_n(\xi)\,d\mu(\xi),
\label{5}
\end{equation}
where the integral converges absolutely.
\end{theorem}
\Proof
If $I$ is bounded then the integral in \eqref{5} converges absolutely
by continuity of $f$. If $I$ is unbounded then, for fixed
$\de\in[0,\ep]$,
we have for some $r\ge0$ that
$f(x+\de\xi)=O(|\xi|^r)$ as $\xi\to\pm\iy$ on $\supp(\mu)$. So also in
that case the integral in \eqref{5} converges absolutely.

By substitution of \eqref{57}, by orthogonality and by \eqref{2}
we have:
\begin{align*}
\frac1{\delta^n}\,\int_\RR f(x+\de\xi)\,p_n(\xi)\,d\mu(\xi)
&=\sum_{k=0}^n\frac{f^{(k)}(x)\,\delta^{k-n}}{k!}\,
\int_\RR \xi^k\,p_n(\xi)\,d\mu(\xi)
+\int_\RR F_{x,n}(\de\xi)\,\xi^n\,p_n(\xi)\,d\mu(\xi)\\
&=\frac{f^{(n)}(x)}{n!}\,\int_\RR \xi^n\,p_n(\xi)\,d\mu(\xi)
+\int_\RR F_{x,n}(\de\xi)\,\xi^n\,p_n(\xi)\,d\mu(\xi)\\
&=\frac{h_n}{k_n n!}\,f^{(n)}(x)
+\int_\RR F_{x,n}(\de\xi)\,\xi^n\,p_n(\xi)\,d\mu(\xi).
\end{align*}
Thus the theorem will be proved if we can show that
\begin{equation}
\lim_{\de\downarrow0}\int_\RR F_{x,n}(\de\xi)\,\xi^n\,p_n(\xi)\,d\mu(\xi)=0.
\label{55}
\end{equation}
By the second part of Proposition \ref{54} we have the estimate
$|F_{x,n}(h)|\le C(1+|h|)^r$ ($h\in I-x$) for some $C>0$, $r\ge0$.
Hence, for $\de\in[0,\ep]$ and $\xi\in\supp(\mu)$ we have the estimate
$|F_{x,n}(\de\xi)|\le C(1+\ep|\xi|)^r$. Thus, the dominated convergence
theorem can be applied to the \LHS\ of \eqref{55}. Then, again by 
Proposition \ref{54}, it follows that \eqref{55} is true.
\qed
\bPP
Note the following special cases of \eqref{5}.
\begin{itemize}
\item
{\bf Gram polynomials} $x\mapsto t_n(x,N)$
(use \eqref{12}):
\begin{equation}
f^{(n)}(x)=\frac{(2n+1)!\,(N-n-1)!}{n!\,(N+n)!}\,
\lim_{\delta\downarrow0}\frac1{\delta^n}\,
\sum_{\xi=0}^{N-1} f(x+\de\xi)\,t_n(\xi,N).
\label{152}
\end{equation}
\item
{\bf Centered Gram polynomials} $x\mapsto t_n(x+N,2N+1)$ on $2N+1$
lattice points (use \eqref{68}):
\begin{equation}
f^{(n)}(x)=\frac{(2n+1)!\,(2N-n)!}{n!\,(2N+n+1)!}\,
\lim_{\delta\downarrow0}\frac1{\delta^n}\,
\sum_{\xi=-N}^N f(x+\de\xi)\,t_n(\xi+N,2N+1).
\label{74}
\end{equation}
In particular, for $n=1$:
\begin{equation}
f'(x)=\frac3{2N(N+\thalf)(N+1)}\,
\lim_{\delta\downarrow0}\frac1\delta\,
\sum_{\xi=-N}^N f(x+\de\xi)\,\xi.
\label{70}
\end{equation}
Analogues of \eqref{74} and \eqref{70} might also be given for
the centered Gram polynomials $x\mapsto t_n(x-N+\thalf,2N)$ on
$2N$ lattice points.
\item
{\bf Legendre polynomials} $P_n$ (use \eqref{13}):
\begin{equation}
f^{(n)}(x)=\frac{(2n+1)!}{2^{n+1}  n!}\,
\lim_{\delta\downarrow0}\frac1{\delta^n}\,
\int_{-1}^1 f(x+\de\xi)\,P_n(\xi)\,d\xi.
\label{73}
\end{equation}
In particular, for $n=1$:
\begin{equation}
f'(x)=\frac32\,\lim_{\delta\downarrow0}\frac1\delta\,
\int_{-1}^1 f(x+\de\xi)\,\xi\,d\xi.
\label{69}
\end{equation}
\item
{\bf Shifted Legendre polynomials} $x\mapsto P_n(2x-1)$
(use \eqref{59}):
\begin{equation}
f^{(n)}(x)=\frac{(2n+1)!}{n!}\,\lim_{\delta\downarrow0}\frac1{\delta^n}\,
\int_0^1 f(x+\de\xi)\,P_n(2\xi-1)\,d\xi.
\label{60}
\end{equation}
\end{itemize}
\subsection{Cioranescu's 1938 paper}
A variant of Theorem \ref{16} was first stated and proved by
Cioranescu \cite[formula $({\rm M}')$]{31} in 1938 for the case that
$d\mu(x)=w(x)\,dx$ is
absolutely continuous with bounded support within an interval $[a,b]$.
He showed for $f\in C^n([a,b])$ that there exists $\eta\in(a,b)$
such that
\begin{equation}
n!\,\frac{\int_a^b f(y)\,p_n(y)\,w(y)\,dy}
{\int_a^b y^n\,p_n(y)\,w(y)\,dy}=f^{(n)}(\eta).
\label{58}
\end{equation}
Then he took limits for $b\downarrow a$ in the \LHS\ of \eqref{58}
(see \cite[formula (9)]{31}) with $p_n$ remaining an orthogonal polynomial
on the shrinking interval $[a,b]$ with respect to the weight function $w$
restricted to $[a,b]$. The limit on the \RHS\ of \eqref{58} then becomes
$f^{(n)}(a)$. In general, this limit formula for $f^{(n)}(a)$ will not
be contained in \eqref{5} since the weight function (after rescaling it
to a fixed interval) will not remain the same during the limit process.
But Cioranescu's limit result in the case of shifted Legendre polynomials
is the same as \eqref{60}. The case $n=1$ of \eqref{60} is explicitly
mentioned by Cioranescu (see \cite[formula $(9')$]{31}).
\subsection{Substitution of the Rodrigues formula}
For classical orthogonal polynomials $p_n$ and for $f\in C^n(I)$ with
$f^{(n)}$ of polynomial growth on $I$ we can also prove
\eqref{5} by substituting $d\mu(x)=w(x)$ and the Rodrigues formula
\eqref{103} together with \eqref{104}, and by
performing integration by parts $n$ times:
\begin{equation}
\frac{k_n n!}{h_n \delta^n}\,
\int_I f(x+\de\xi)\,p_n(\xi)\,w(\xi)\,d\xi
=\frac{\int_I f^{(n)}(x+\de\xi)\,X(\xi)^n\,w(\xi)\,d\xi}
{\int_I X(\xi)^n\,w(\xi)\,d\xi}\,\to f^{(n)}(x)\quad
\mbox{as $\de\downarrow0$.}
\label{108}
\end{equation}
In the Legendre case $w(x)=1$ this was already observed by
Cioranescu \cite[p.296]{31}.
\subsection{Haslam-Jones' 1953 paper}
Next Theorem \ref{16}, for the case that $\mu$ has bounded support,
was observed (with proof omitted as being easy)
in 1953 by Haslam-Jones \cite[p.192]{26}, who was apparently not aware of
Cioranescu's result.
In fact, in his formulation the measure
$\mu$ only has to be real, not necessarily positive.
Furthermore, $f$ only has to be continuous with an $n$-th Peano derivative
at $x$ (see \eqref{76}). Note that our proof of Theorem \ref{16}
can be used without essential changes under the weaker hypotheses of
Haslam-Jones.

In fact, the assumptions in \cite{26} are still weaker. Haslam-Jones
assumes, for given $n>0$, a real, not necessarily positive measure $\nu$
with bounded support (or equivalently a function $\nu$ of bounded variation)
on a finite interval $J$ such that
$\int_J x^k\,d\nu(x)=0$ for $k=0,1,\ldots,n-1$
and $\int_J x^n\,d\nu(x)=\ka\ne0$. Then for a function $f$ which
is continuous on a neighbourhood of $x$ and has $n$-th Peano derivative
$c_n$ in $x$ we have
\begin{equation}
c_n=\frac{n!}{\ka}\,\lim_{\delta\downarrow0}\frac1{\delta^n}\,
\int_J f(x+\de\xi)\,d\nu(\xi).
\label{61}
\end{equation}
Again this can be proved as we did for Theorem \ref{16}, without essential
changes.
\subsection{A special case of Haslam-Jones' results}
We will consider here a special case of \eqref{61} which is essentially
different from Theorem \ref{16}. Let $\{p_m\}$ be a system of
orthogonal polynomials on $[-1,1]$ with
respect to a positive Borel measure~$\mu$. Let $K_m(x,y)$ be the
corresponding Christoffel-Darboux kernel given by \eqref{31}, \eqref{62}.
Fix $n$ and define the measure $\nu$ in \eqref{61} by
\begin{equation}
\int_{-1}^1 f(\xi)\,d\nu(\xi):=f(1)
-\int_{-1}^1 f(\xi)\,K_{n-1}(\xi,1)\,d\mu(\xi).
\label{65}
\end{equation}
Indeed, by the reproducing kernel property the \RHS\
of \eqref{65} equals 0 if $f$ is
a polynomial of degree $<n$, while for $f(\xi):=\xi^n$ the
\RHS\ of \eqref{65} becomes
\[
-\int_{-1}^1 \xi^{n-1}(\xi-1)\,K_{n-1}(\xi,1)\,d\mu(\xi)=
\frac{k_{n-1}p_n(1)}{k_n h_{n-1}}\,
\int_{-1}^1 \xi^{n-1}\,p_{n-1}(\xi)\,d\mu(\xi)
=\frac{p_n(1)}{k_n}\ne0.
\]
Thus for this case \eqref{61} becomes
\begin{equation}
c_n=\frac{n!\,k_n}{p_n(1)}\,\lim_{\delta\downarrow0}
\frac1{\delta^n}\,\left(f(\de)
-\int_{-1}^1 f(x+\de\xi)\,K_{n-1}(\xi,1)\,d\mu(\xi)\right).
\label{63}
\end{equation}
In particular, take $d\mu(\xi):=d\xi$ on $[-1,1]$.
Then substitute \eqref{67}, by which
\eqref{63} takes the form
\begin{equation}
c_n=2^n(\thalf)_n\,\lim_{\delta\downarrow0}
\frac1{\delta^n}\,\left(f(\de)
-\thalf\int_{-1}^1 f(x+\de\xi)\,(P_n'(\xi)+P_{n-1}'(\xi))\,d\xi\right).
\label{64}
\end{equation}
For this case Haslam-Jones \cite{26} showed that if the limit
on the right of \eqref{64} exists then the $n$-th Peano derivative of
$f$ at $x$ exists and it equals $c_n$ given by \eqref{64}.
A different proof of this result was given by Gordon \cite{35}.
\subsection{Connection with Jacobi type orthogonal polynomials}
For a larger family of examples than \eqref{64} consider formula
\eqref{65} with
$d\mu(x):=(1-x)^\al(1+x)^{\be+1}\,dx$ ($\al,\be>-1$). Then
$p_m(x)=P_m^{(\al,\be+1)}(x)$, a {\em Jacobi polynomial}
(see \cite[Ch.~4]{1}). From \cite[(4.5.3)]{1} we obtain that
\[
{\bf K}_{n-1}(x,1)
=\frac{\Ga(n+\al+\be+2)}{2^{\al+\be+2}\Ga(\al+1)\Ga(n+\be+1)}\,
P_{n-1}^{(\al+1,\be+1)}(x).
\]
Then the vanishing of the \RHS\ of \eqref{65} for polynomials $f$ of
degree $<n$ can be written more explicitly as
\begin{multline*}
\frac{\Ga(\al+\be+2)}{2^{\al+\be+1}\Ga(\al+1)\Ga(\be+1)}\,
\int_{-1}^1 f(x)\,
\frac{(1+x)\,P_{n-1}^{(\al+1,\be+1)}(x)}{2P_{n-1}^{(\al+1,\be+1)}(1)}\,
(1-x)^\al(1+x)^\be\,dx\\
-\frac{(\be+1)_n (n-1)!}{(\al+\be+2)_n (\al+2)_{n-1}}\,f(1)=0.
\end{multline*}
Hence, for fixed $n$, the polynomial $(1+x) P_{n-1}^{(\al+1,\be+1)}(x)$
is the $n$-th degree orthogonal polynomial with respect to a measure
on $[-1,1]$ consisting of the weight function $(1-x)^\al (1+x)^\be$
and a negative constant times a delta weight at $x=1$.
On comparing with \cite[Theorem 3.1]{36} (extended to negative multiples
of delta weights by analytic continuation) we can identify this
orthogonal polynomial with a {\em Jacobi type polynomial}\,:
\[
P_n^{(\al,\be;0,N)}(x)=c(1+x)P_{n-1}^{(\al+1,\be+1)}(x),
\]
where
\[
N=-\,\frac{(\be+1)_n (n-1)!}{(\al+\be+2)_n (\al+2)_{n-1}}\quad
\mbox{and}\quad
c=\frac{P_n^{(\al,\be;0,N)}(1)}{2P_{n-1}^{(\al+1,\be+1)}(1)}=
\frac{\al+1}{2n}\,.
\]
This corresponds correctly with \cite[(2.1)]{36}, which simplifies for $M=0$
and the above choice of $N$ to
\[
P_n^{(\al,\be;0,N)}(x)=\frac{\al+1}{n(n+\al+\be+1)}\,(1+x)\,
\frac d{dx} P_n^{(\al,\be)}(x).
\]
The above formulas extend by continuity to the case $\be=-1$, which
occurs if $d\mu(x)=dx$ (the case considered in \eqref{64}). 
\subsection{Lanczos' 1956 work}
In a book published in 1956 Lanczos \cite[(5-9.1)]{19}, apparently
unaware of the earlier work
by Cioranescu \cite{31} and Haslam-Jones \cite{26}, rediscovered
formula \eqref{69}.
He called this {\em differentiation by integration}.
His work got quite a lot of citations, see for instance
\cite{20}, \cite{37}, \cite{55}, \cite{24}, \cite{21}, \cite{22}, \cite{25}.
The name {\em Lanczos derivative},
notated as $f_L'(x)$, became common for a value obtained from
the \RHS\ of \eqref{69}.
In \cite{24} and \cite{21} also
\eqref{73} (the Legendre case for general $n$) was rediscovered.

As an important new aspect Lanczos observed that
\[
\frac3{2\de}\,\int_{-1}^1 f(x+\de\xi)\,\xi\,d\xi
=
\frac{\int_{-1}^1 f(x+\de\xi)\,\xi\,d\xi}
{\de\int_{-1}^1 \xi^2\,d\xi},
\]
as an approximation of $f'(x)$, is the limit as $N\to\iy$
of the quotient of Riemann sums
\[
\frac{N^{-1}\sum_{\xi=-N}^N f(x+\de\xi/N)\,\xi/N}
{\de N^{-1}\sum_{\xi=-N}^N (\xi/N)^2}=
\frac3{2N(N+\thalf)(N+1)}\,\frac N\de\sum_{\xi=-N}^N f(x+N^{-1}\de\xi)\,\xi.
\]
We have seen this last expression already in \eqref{70} (the special
case of \eqref{5} with a centered Gram polynomial of degree 1).
The expression approximates
$f'(x)$ for $N^{-1}\de$ small.
In fact, \eqref{70} was the starting point of Lanczos, see \cite[(5-8.4)]{19}.
\subsection{Interpretation by least-square approximation}
Lanczos \cite[(5-8.4)]{19} arrived at \eqref{70} by a least-square
minimization (a technique invented by Gauss and Legendre, see for instance
\cite{59}). For this compare the function $\xi\mapsto f(x+\de\xi)$
with a linear function
$g(\xi):=a_0+a_1\xi$ such that the squared distance
\[
S(a_0,a_1):=\sum_{\xi=-N}^N \big(f(x+\de\xi)-g(\xi)\big)^2
=\sum_{\eta=-\de N,-\de(N-1),\ldots,\de N}
\big(f(x+\eta)-(a_0+\de^{-1}a_1\eta)\big)^2
\]
is minimal. Then the slope $\de^{-1}a_1$ of the straight line
$\eta\mapsto a_0+\de^{-1}a_1\eta$ minimizing the
distance will approximate $f'(x)$ for small $\de$.
The minimum is achieved for a unique $(a_0,a_1)$, where
one finds $a_1$ in this simple case already by solving
$\frac\pa{\pa a_1}S(a_0,a_1)=0$. Thus Lanczos obtained
\[
\de^{-1}a_1=\frac{\sum_{\xi=-N}^N f(x+\de\xi)\xi}
{\de\sum_{\xi=-N}^N \xi^2}\,=
\frac3{2\de N(N+\thalf)(N+1)}\,
\sum_{\xi=-N}^N f(x+\de\xi)\,\xi,
\]
and he thus arrived at \eqref{70}.

We can interpret \eqref{5} as a more general least-square approximation.
By the assumptions on $f$ in the Theorem \ref{16}
the function $\xi\mapsto f(x+\xi\delta)$ is in $L^2(\RR,\mu)$
for each $\delta\in[0,\ep]$. Let $\xi\mapsto P_{n,x,\delta}[f](\xi)$
be the
polynomial of degree $\le n$ which is on minimal distance from the
function
$\xi\mapsto f(x+\xi\delta)$ in the Hilbert space $L^2(\RR,\mu)$. Then
\begin{equation}
P_{n,x,\delta}[f](\eta)=\sum_{k=0}^n\frac1{h_k}
\Big(\int_\RR f(x+\delta\xi)\,p_k(\xi)\,d\mu(\xi)\Big)\,p_k(\eta).
\label{24}
\end{equation}
Then
\[
\de^{-n}(P_{n,x,\delta}[f])^{(n)}(0)=\frac{p_n^{(n)}(0)}{\de^nh_n}\,
\int_\RR f(x+\delta\xi)\,p_n(\xi)\,d\mu(\xi)
=\frac{k_n n!}{\de^nh_n}\,\int_\RR f(x+\delta\xi)\,p_n(\xi)\,d\mu(\xi)
\]
approximates $f^{(n)}(x)$ as $\de\downarrow0$. Thus we arrive
at \eqref{5}.

Also observe that clearly
\[
\int_\RR f(x+\delta\xi)\,p_n(\xi)\,d\mu(\xi)=
\int_\RR P_{n,x,\delta}[f](\xi)\,p_n(\xi)\,d\mu(\xi).
\]
Hence, we can rewrite \eqref{5} as
\begin{equation}
f^{(n)}(x)=\frac{k_n n!}{h_n}\,\lim_{\delta\downarrow0}\frac1{\delta^n}\,
\int_\RR P_{n,x,\delta}[f](\xi)\,p_n(\xi)\,d\mu(\xi).
\label{6}
\end{equation}

For the Legendre case \eqref{73} this interpretation by
least-square approximation was given in~\cite{21}.
But earlier, in 1990, Kopel \& Schramm \cite{38}, apparently unaware of
any predecessors, arrived at the case $n=1$ of \eqref{60} while they were
guided by least-square approximation.
\subsection{Even orthogonality measures}
We can refine Theorem \ref{16} if we assume that the orthogonality
measure $\mu$ considered there is {\em even}, i.e., that $d\mu(-x)=d\mu(x)$.
Then the corresponding orthogonal polynomials $p_n$ are even or odd
according to whether $n$ is even or odd, respectively. The simplest
example is given by the Legendre polynomials. Now also modify the
assumptions about $f^{(n)}(x)$ in Theorem \ref{16}. Only assume that
the right $n$-th derivative $f_+^{(n)}(x)$ and left
$n$-th derivative $f_-^{(n)}(x)$ exist. Then by an easy adaptation
of the proof of Theorem \ref{16} we get \eqref{16} with on the \LHS\
the symmetric $n$-th derivative of $f$ at $x$:
\begin{equation}
\thalf\big(f_+^{(n)}(x)+f_-^{(n)}(x)\big)
=\frac{k_n n!}{h_n}\,\lim_{\delta\downarrow0}\frac1{\delta^n}\,
\int_\RR f(x+\de\xi)\,p_n(\xi)\,d\mu(\xi).
\label{71}
\end{equation}
The special case of this result for Legendre polynomials (see
\eqref{73}, \eqref{69}) was observed in \cite[Proposition 1]{20}
for $n=1$ and in \cite[Theorem 2]{21} for general $n$.
The special case of \eqref{71}
for centered Gram polynomials (see \eqref{74})
was observed in \cite[Theorem 1.1]{22}.

Moreover, in \cite[pp.~370--371]{21} and
\cite[\S4]{25} examples were given for the Legendre case with $n=1$,
where the limit on the \RHS\ of \eqref{69} exists, but the left and right
derivative of $f$ at $x$ do not exist.
Earlier, in \cite[pp.~231--232]{38} an example was given where the
limit of the \RHS\ of \eqref{60} for $n=1$ exists, while the right
derivative of $f$ at $x$ does not exist.

Consider the proof of Theorem \ref{16} if we know that $f^{(k)}(x)$ exists
for $k$ up to some $m>n$. Then
\begin{multline*}
\frac1{\delta^n}\,\int_\RR f(x+\de\xi)\,p_n(\xi)\,d\mu(\xi)
=\frac{h_n}{k_n n!}\,f^{(n)}(x)
+\sum_{k=m+1}^n\frac{f^{(k)}(x)\,\delta^{k-n}}{k!}\,
\int_\RR \xi^k\,p_n(\xi)\,d\mu(\xi)\\
+\de^{m-n}\int_\RR F_{x,m}(\de\xi)\,\xi^m\,p_n(\xi)\,d\mu(\xi)
=\frac{h_n}{k_n n!}\,f^{(n)}(x)+\FSO(\de)\quad
\mbox{as $\de\downarrow0$.}
\end{multline*}
We can say more if moreover $\mu$ is an even measure.
Then $\int_\RR \xi^{n+1}\,p_n(\xi)\,d\mu(\xi)=0$ and thus we have,
as $\de\downarrow0$,
\begin{equation}
\frac1{\delta^n}\,\int_\RR f(x+\de\xi)\,p_n(\xi)\,d\mu(\xi)
=\frac{h_n}{k_n n!}\,f^{(n)}(x)+
\begin{cases}
o(\de)&\mbox{if $f^{(n+1)}(x)$ exists,}\\
\FSO(\de^2)&\mbox{if $f^{(n+2)}(x)$ exists.}
\end{cases}
\label{72}
\end{equation}
The Legendre case of \eqref{72}
was observed for $n=1$ in \cite[(5-9.3)]{19} and for general $n$ in
\cite[(4)]{24}.
\subsection{Generalized Taylor series}
Rewrite \eqref{5} as
\begin{equation}
f^{(n)}(x)=\frac{k_n n!}{h_n}\,\lim_{\delta\downarrow0}\frac1{\delta^n}\,
\int_\RR f(x+\xi)\,p_n(\de^{-1}\xi)\,d\mu_\de(\xi).
\label{105}
\end{equation}
Here $\mu_\de(E):=\mu(\de^{-1}E)$. Then the polynomials 
$x\mapsto p_n(\de^{-1}x)$ are orthogonal with respect to the measure
$\mu_\de$. The formal Taylor series
\begin{equation}
\sum_{n=0}^\iy \frac{f^{(n)}(x)}{n!}\,\eta^n
\label{107}
\end{equation}
can accordingly be seen as a termwise limit of the formal generalized
Fourier series
\begin{multline}
\sum_{n=0}^\iy \frac1{h_n}
\left(\int_\RR f(x+\xi)\,p_n(\de^{-1}\xi)\,d\mu_\de(\xi)\right)p_n(\eta)\\
=\sum_{n=0}^\iy\frac1{n!}
\left(\frac{k_n\,n!}{\de^n\,h_n}\,
\int_\RR f(x+\xi)\,p_n(\de^{-1}\xi)\,d\mu_\de(\xi)\right)
\frac{\de^n\,p_n(\de^{-1}\eta)}{k_n}\,.
\label{106}
\end{multline}
Indeed, use \eqref{105} and the limit
$\de^n\,p_n(\de^{-1}\eta)/k_n\to \eta^n$ as $\de\downarrow0$.

While for a big class of orthogonal polynomials
and for $f$ moderately smooth, the series \eqref{106}
converges with sum $f(x+\eta)$ because of equiconvergence theorems
given in Szeg{\H{o}} \cite[Ch.~9 and 13]{1}, we will need analyticity
of $f$ in a neighbourhood of $x$ for convergence of \eqref{107} to
$f(x+\eta)$. In the case of Jacobi polynomials and for $f$ analytic
on a neighbourhood of $x$ we can use Szeg{\H{o}} \cite[Theorem 9.1.1]{1}.
Then an open disk around $x$ of radius less than the convergence radius
of \eqref{107} is in the interior of the ellipse of convergence
\eqref{106} for $\de$ small enough. This was discussed for Legendre
polynomials by Fishback \cite{56}.

A different limit from orthogonal polynomials to monomials is discussed
by Askey \& Haimo~\cite{57}. This involves Gegenbauer polynomials,
which we write as Jacobi polynomials
$P_n^{(\al,\al)}(x)=k_n^{(\al,\al)} x^n+\cdots\;$:
\begin{equation}
\lim_{\al\to\iy}\frac{P_n^{(\al,\al)}(x)}{k_n^{(\al,\al)}}=x^n.
\label{109}
\end{equation}
Consider now the formal expansion of $f(x+\eta)$ for
$\eta\in[-1,1]$ in terms of the polynomials $P_n^{(\al,\al)}(\eta)$:
\begin{align*}
&\sum_{n=0}^\iy \frac1{h_n^{(\al,\al)}}
\left(\int_{-1}^1 f(x+\xi)\,P_n^{(\al,\al)}(\xi)\,(1-\xi^2)^\al\,d\xi\right)
P_n^{(\al,\al)}(\eta)\\
&=\sum_{n=0}^\iy\frac1{n!}
\left(\frac{k_n^{(\al,\al)}\,n!}{h_n^{(\al,\al)}}\,
\int_{-1}^1 f(x+\xi)\,P_n^{(\al,\al)}(\xi)\,(1-\xi^2)^\al\,d\xi\right)
\frac{P_n^{(\al,\al)}(\eta)}{k_n^{(\al,\al)}}\\
&=\sum_{n=0}^\iy\frac1{n!}\,
\frac{\int_{-1}^1 f^{(n)}(x+\xi)\,(1-\xi^2)^{n+\al}\,d\xi}
{\int_{-1}^1 (1-\xi^2)^{n+\al}\,d\xi}\,
\frac{P_n^{(\al,\al)}(\eta)}{k_n^{(\al,\al)}}\,,
\end{align*}
where we used the identity for $\de=1$ in \eqref{108}.
The last form of the above formal series tends termwise to the formal
Taylor series \eqref{107}. This is seen from \eqref{109} and the fact
that the measure
$(1-\xi^2)^{n+\al}\,d\xi/\int_{-1}^1 (1-\xi^2)^{n+\al}\,d\xi$
tends to the delta measure as $\al\to\iy$ (see \cite[p.301]{57}).
As observed in \cite[p.303]{57}, the function $f$ has to be increasingly
smooth as $\al$ grows in order to have convergence in the expansion
of $f(x+\eta)$ in terms of $P_n^{(\al,\al)}(\eta)$.
\subsection{Connection with the continuous wavelet transform}
\label{17}  
The {\em continuous wavelet transform} $\Phi_g$
(see for instance \cite{12}, \cite{13}) is defined by
\begin{equation}
(\Phi_g f)(a,b):=|a|^{-\half}\,
\int_\RR f(t)\,\overline{g\big(a^{-1}(t-b)\big)}\,dt
\qquad(f\in L^2(\RR),\;a,b\in\RR,\;a\ne0).
\label{20}
\end{equation}   
Here we will take the {\em wavelet} $g$ as a nonzero function in
$(L^1\cap L^2)(\RR)$ such that $\int_\RR g(t)\,dt=0$.

For orthogonal polynomials $p_n(x)$ let the orthogonality measure have the form
$d\mu(x)=w(x)\,dx$ (essentially item 1 in \S\ref{18})
with $w(x)\ge0$ for $x\in\RR$ and with the functions $x\mapsto x^n\,w(x)$
in $(L^1\cap L^2)(\RR)$ for all $n=1,2,\ldots\;$.
Put
\begin{equation}
g_n(x):=p_n(x)\,w(x).
\label{21}
\end{equation}
Then, for $n=1,2,\ldots\;$ the functions $g_n$ are in $(L^1\cap
L^2)(\RR)$ and satisfy $\int_\RR g_n(t)\,dt=0$.
Now consider the continuous wavelet transform \eqref{20} for $g$ equal
to such $g_n$ and compare with \eqref{5}.
Then
\begin{equation}
\int_\RR f(x+\xi\delta)\,p_n(\xi)\,w(\xi)\,d\xi=
\delta^{-1}\,\int_\RR f(t)\,g_n\big(\delta^{-1}(t-x)\big)\,dt=
\delta^{-\half}\,(\Phi_{g_n}f)(\delta,x).
\end{equation}
Hence, \eqref{5} can now be written as
\begin{equation}
f^{(n)}(x)
=\frac{k_n n!}{h_n}\,\lim_{\delta\downarrow0}\,\frac1{\delta^{n-\half}}\,
(\Phi_{g_n}f)(\delta,x).
\end{equation}

A similar observation about the continuous wavelet transform
approximating the $n$-th derivative was made by Rieder \cite[(5)]{11}
in the case of \eqref{20} with $g$ having its first $n$ moments equal
to zero.
In fact, he recovers formula \eqref{61}, first obtained by
Haslam-Jones \cite{26}, for $\nu$ absolutely continuous,
$f$, $c_n$ being the $n$-th distributional derivative of $f$,
and the limit taken in a suitable Sobolev norm (see also
\cite[Theorem 2.3]{11}). An example of a wavelet $g$ having its first
$n$ moments equal to zero is Daubechies' wavelet ${}_n\psi$
of compact support, see \cite[p.984]{52}.

The continuous wavelet transform $\Phi_{g_n}$ with $g_n$ given by
\eqref{21} and $w(x)$ being a weight function for one of the classical
orthogonal polynomials (Jacobi, Laguerre, Hermite) was considered
by Moncayo \& Y\'a\~nez \cite{15}.
\subsection{A special case: the $n$-th order finite difference as approximation of
the $n$-th derivative}
For $N=n+1$ we can see that \eqref{152} specializes as an $n$-th order finite
difference approximating the $n$-th derivative. Indeed, write \eqref{152} as
\begin{equation*}
f^{(n)}(x)=\lim_{\de\downarrow0}(D_{n,\de,N}f)(x),
\end{equation*}
where (see also \eqref{11} and \cite[(9.5.1)]{5})
\begin{align*}
(D_{n,\de,N}f)(x)&:=
\frac{(2n+1)!\,(N-n-1)!}{n!\,(N+n)!}\,
\frac1{\delta^n}\,
\sum_{\xi=0}^{N-1} f(x+\de\xi)\,t_n(\xi,N)\\
&\;=
\frac{(2n+1)!\,(N-1)!\,(-1)^n}{n!\,(N+n)!\,\de^n}\,
\sum_{\xi=0}^{N-1} f(x+\de\xi)\,\hyp32{-n,n+1,-\xi}{1,-N+1}1.
\end{align*}
Now put $N:=n+1$ and use that for $\xi=0,1,\ldots,n$ we have by
\cite[(15.2.4)]{53}:
\[
\hyp32{-n,n+1,-\xi}{1,-n}1=\sum_{k=0}^\xi \frac{(-\xi)_k\,(n+1)_k}{k!\,k!}
=\hyp21{-\xi,n+1}11=\frac{(-n)_\xi}{(1)_\xi}=(-1)^\xi\binom n\xi.
\]
Hence
\[
(D_{n,\de,n+1}f)(x)=\de^{-n}\sum_{\xi=0}^n (-1)^{n-\xi}\binom n\xi f(x+\de\xi)
=(\De_\de^n f)(x),
\]
where
\[
(\De_\de f)(x):=\frac{f(x+\de)-f(x)}\de\,.
\]
\section{Filters for higher derivatives}
\label{137}
The following theorem is a multi-term variant of Theorem \ref{16}.
\begin{theorem}
\label{37}
Let $\{p_k\}_{k=0,1,2,\ldots}$ be a system of orthogonal polynomials with
respect to the orthogonality measure $\mu$. Let $m,n$ be integers such that
$0\le m\le n$. Let $x\in\RR$.
Let $I$ be a closed interval such that,
for some $\ep>0$, $x+\delta\xi\in I$
if $0\le \delta\le\ep$ and $\xi\in\supp(\mu)$.
Let $f$ be a continuous function on $I$ such that its derivatives of order
$1,2,\ldots,n$ at $x$ exist. In addition, if
$I$ is unbounded, assume that $f$ is of at most polynomial
growth on $I$. Then
\begin{equation}
\frac1{\delta^m}\,\sum_{j=m}^n\frac1{h_j}
\left(\int_\RR f(x+\de\xi)\,p_j(\xi)\,d\mu(\xi)\right)p_j^{(m)}(0)=
f^{(m)}(x)+o(\de^{n-m})\quad\mbox{as $\de\downarrow0$}.
\label{29}
\end{equation}
If $f\in C^n(I)$ and $f^{(n)}$ is of polynomial growth on $I$ then
\eqref{29} holds uniformly for $x$ in compact subsets of $I$.
\end{theorem}
Note that \eqref{29} turns down to \eqref{5} if $m=n$.
\bLP
\Proof
With the notation \eqref{31} we can rewrite the \LHS\ of \eqref{29} as
\begin{equation}
\de^{-m}\,\left(\frac\pa{\pa\eta}\right)^m\,
\int_{\RR} f(x+\de\xi)\,{\bf K}_n(\xi,\eta)\,d\mu(\xi)\Big|_{\eta=0}.
\label{38}
\end{equation}
By substitution of \eqref{57} and by application of \eqref{30}
the expression \eqref{38} is equal to
\begin{multline*}
\sum_{l=0}^{n-m}\frac{f^{(m+l)}(x)}{l!}\,\de^l\eta^l\Big|_{\eta=0}+
\de^{n-m}\,\left(\frac\pa{\pa\eta}\right)^m\,
\int_{\RR} \xi^n\,F_{x,n}(\de\xi)\,
{\bf K}_n(\xi,\eta)\,d\mu(\xi)\Big|_{\eta=0}\\
=f^{(m)}(x)+\de^{n-m}\int_{\RR} \xi^n\,F_{x,n}(\de\xi)\,
\left(\frac\pa{\pa\eta}\right)^m\,K_n(\xi,\eta)\Big|_{\eta=0}\,d\mu(\xi).
\end{multline*}
Now use the same argument as at the end of the proof of Theorem \ref{16},
calling Proposition \ref{54} and using dominated convergence,
in order to show that
\[
\qquad\qquad\qquad\qquad\quad
\lim_{\de\downarrow0}\int_{\RR} \xi^n\,F_{x,n}(\de\xi)\,
\left(\frac\pa{\pa\eta}\right)^m\,
{\bf K}_n(\xi,\eta)\Big|_{\eta=0}\,d\mu(\xi)=0.
\]
For the proof of the last statement also use Proposition \ref{81}.
\qed
\begin{remark}
In view of \eqref{24} we can rewrite \eqref{29} as
\begin{equation}
f^{(m)}(x)=\frac1{\de^m}\,
\frac{d^m}{d\eta^m}\,P_{n,x,\de}[f](\eta)\Big|_{\eta=0}+o(\de^{n-m})
\quad\mbox{as $\de\downarrow0$}\qquad(m=0,1,\ldots,n).
\label{25}
\end{equation}
If $f$ is a polynomial of degree $\le n$ then \eqref{25} holds exactly
without the term $o(\de^{n-m})$ because
\[
P_{n,x,\de}[f](\eta)=f(x+\de\eta).
\]
\end{remark}

For $f$ having derivatives at $x$ up to order $n+1$
we can refine \eqref{29} as follows.
\begin{proposition}
\label{41}
Keep the assumptions of Theorem \ref{37}. Moreover assume that
$f^{(n+1)}(x)$ exists. Then
\begin{multline}
\frac1{\delta^m}\,\sum_{j=m}^n\frac1{h_j}
\left(\int_\RR f(x+\de\xi)\,p_j(\xi)\,d\mu(\xi)\right)p_j^{(m)}(0)\\
=
f^{(m)}(x)-\frac{p_{n+1}^{(m)}(0)\,f^{(n+1)}(x)}{k_{n+1}\,(n+1)!}\,\de^{n-m+1}
+o(\de^{n-m+1})\quad\mbox{as $\de\downarrow0$}.
\label{39}
\end{multline}
If $f\in C^{n+1}(I)$ and $f^{(n+1)}$ is of polynomial growth on $I$ then
\eqref{39} holds uniformly for $x$ in compact subsets of $I$.
\end{proposition}
\Proof
Write the \LHS\ of \eqref{39} as \eqref{38}.
Then substitute \eqref{57} with $n$ replaced by $n+1$
and apply \eqref{30}. Then the expression \eqref{38} becomes
\begin{multline*}
f^{(m)}(x)+
\frac{f^{(n+1)}(x)}{(n+1)!}\,\de^{n-m+1}\,
\left(\frac\pa{\pa\eta}\right)^m\,
\int_\RR\xi^{n+1}\,{\bf K}_n(\xi,\eta)\,d\mu(\xi)\,\Big|_{\eta=0}\\
+\de^{n-m+1}\,\int_{\RR} \xi^{n+1}\,F_{x,n+1}(\de\xi)\,
\left(\frac\pa{\pa\eta}\right)^m\,
{\bf K}_n(\xi,\eta)\Big|_{\eta=0}\,d\mu(\xi).
\end{multline*}
By a similar argument as in the proof of Theorem \ref{37} we see that
the last term equals $o(\de^{n-m+1})$ as $\de\downarrow0$.
By application of \eqref{30} the second term becomes
\begin{multline*}
\frac{f^{(n+1)}(x)}{(n+1)!}\,\de^{n-m+1}\,
\left(\frac\pa{\pa\eta}\right)^m\,\int_\RR\xi^{n+1}\,
\left({\bf K}_{n+1}(\xi,\eta)-\frac{p_{n+1}(\xi)\,
p_{n+1}(\eta)}{h_{n+1}}\right)
d\mu(\xi)\,\Big|_{\eta=0}\\
=-\frac{p_{n+1}^{(m)}(0)}{h_{n+1}\,(n+1)!}\,\de^{n-m+1}\,
\int_\RR\xi^{n+1}\,p_{n+1}(\xi)\,d\mu(\xi).
\end{multline*}
Then \eqref{39} follows by using \eqref{2}.\qed
\subsection{Even orthogonality measures}
If in Proposition \ref{41} the orthogonality measure $\mu$ is even and if
$n-m$ is odd then $p_j^{(m)}(0)=0$ whenever $j-m$ is odd,
so the sum on the \LHS\
of \eqref{39} then runs over $j=m,m+2,\ldots,n-1$. Moreover, for
$p_{n+1}^{(m)}(0)/k_{n+1}$ occurring on the \RHS\ of \eqref{39}
we then have that
\begin{equation}
(-1)^{(n-m+1)/2}\,\frac{p_{n+1}^{(m)}(0)}{k_{n+1}}>0.
\label{40}
\end{equation}
In order to prove \eqref{40} we can assume without loss of generality
that $k_n>0$. Then, for all $j$ and for $x$ large, $p_{n+1}^{(j)}(x)>0$.
Also note that $p_{n+1}^{(j)}$ is a polynomial of degree $n+1-j$ which
is even or odd according to whether $n+1-j$ is even or odd.
Since $p_{n+1}$ belongs to a family of orthogonal polynomials, it has
$n+1$ simple real zeros. The number of positive zeros is $(n+1)/2$ if
$n+1$ is even and $n/2$ if $n+1$ is odd. Now it follows by induction with
respect to $j$ that $p_{n+1}^{(j)}$ has $(n+1-j)/2$ positive zeros if
$n+1-j$ is even and $(n-j)/2$ if $n+1-j$ is odd. So we arrive also at
this property for $j=m$, and then \eqref{40} readily follows.

Thus, for $\mu$ an even measure and $n-m$ an odd number we can rewrite
\eqref{39} as follows for $\de\downarrow0$.
\begin{multline}
\frac1{\delta^m}\,\sum_{j=m,m+2,\ldots,n-1}
\left(\int_\RR f(x+\de\xi)\,p_j(\xi)\,d\mu(\xi)\right)
\frac{p_j^{(m)}(0)}{h_j}=f^{(m)}(x)
\\-(-1)^{(n-m+1)/2}\,
\frac{|p_{n+1}^{(m)}(0)|\,f^{(n+1)}(x)}{|k_{n+1}|\,(n+1)!}\,
\de^{n-m+1}+\begin{cases}
o(\de^{n-m+1})&\mbox{if $f^{(n+1)}(x)$ exists,}\\
o(\de^{n-m+2})&\mbox{if $f^{(n+2)}(x)$ exists,}\\
\FSO(\de^{n-m+3})&\mbox{if $f^{(n+3)}(x)$ exists.}
\end{cases}\quad
\label{79}
\end{multline}
This is proved by a slight adaptation of the proof of Proposition \ref{41}.
Furthermore, \eqref{79} holds uniformly for $x$ in compact subsets of $I$
if the appropriate derivative of $f$
of order $n+1$, $n+2$ or $n+3$ is continuous
and of polynomial growth on $I$.
\begin{remark}
\label{147}
Note that for fixed $m$ the approximation of the \LHS\ of \eqref{79}
(and earlier formulas \eqref{29}, \eqref{25} and \eqref{39}) to
$f^{(m)}(x)$ becomes better as $n$ increases. However, this observation
disregards the frequency spectrum of the signal $f$ and the effect of noise.
See Remark \ref{148} for a discussion of these aspects.
\end{remark}
\subsection{Filters}
We can consider the \LHS\ of \eqref{29} as a filter (continuous or discrete
depending on the choice of $\mu$) for $m$-th order
differentiation at $x$. In general, a continuous respectively analog
filter sends an input function $f$ to an output function $g$ by
convolution with a fixed real-valued function $\rho$:
\begin{align}
g(y)&=\sum_{x=M}^N f(y-x)\,\rho(x)\qquad(y\in\ZZ),\label{123}\\
g(y)&=\int_{M}^N f(y-x)\,\rho(x)\,dx\qquad(y\in\RR).\label{124}
\end{align}
Here $M$ may be $-\iy$ and $N$ may be $\iy$.
Filters are widely used in electrical engineering, with analog filters
being continuous and digital filters being discrete. There $y$ is usually
the time variable $t$ and instead of $\rho$ one writes $h$, the
{\em unit impulse response}. If $M$ and $N$ are finite in \eqref{123},
one speaks about a {\em finite impulse response} (FIR) filter,
otherwise about an {\em infinite impulse response} (IIR) filter.
The Fourier transform $\phi$ of $\rho$ is called the
{\em characteristic function} or
{\em transfer function} (also denoted by $H$) of the filter:
\begin{align}
\phi(\om)&:=\sum_{x=-M}^N \rho(x)\,e^{-ix\om},
\label{121}\\
\phi(\om)&:=\int_{-M}^N \rho(x)\,e^{-ix\om}\,dx.
\label{122}
\end{align}
Equivalently, $\phi$ equals the quotient $g/f$ of the input function
$g$ and the output function $f$ if $f(y):=e^{i\om y}$.
Standard books of digital filter theory are for instance \cite{58}
and \cite{28}. In \cite[p.306]{58} methods are described for the
design of digital differentiators.
\subsection{The characteristic function}
\label{101}
We continue with the \LHS\ of \eqref{29} considered
as a filter.
Let us assume that $\mu$ is an even measure and that
$n-m$ is odd, so that we can work with \eqref{79}.
We obtain the {\em characteristic  function} $\phi$ of the filter defined
by the \LHS\ of \eqref{79} if we put there
$f(\xi):=e^{i\om \xi}$ with $\om\in\RR$
and take $x:=0$:
\begin{multline}
\phi(\om\de):=\frac1{\delta^m}\,\sum_{j=m,m+2,\ldots,n-1}
\left(\int_\RR e^{i\om\de\xi}\,p_j(\xi)\,d\mu(\xi)\right)
\frac{p_j^{(m)}(0)}{h_j}\\
=(i\om)^m\Big(1-\frac{|p_{n+1}^{(m)}(0)|}{|k_{n+1}|\,(n+1)!}\,
(\de\om)^{n-m+1}+(\de\om)^{n-m+3}\,G(\de\om)\Big),
\label{80}
\end{multline}
with $G$ a bounded function on $\RR$ (for the proof of the
second equality use Proposition \ref{54} with
$f$ bounded).
\begin{remark}
\label{148}
From the last part of \eqref{80} we see that for differentiation
with fixed order $m$ the
first term gives the characteristic function $(i\om)^m$ of the
ideal differentiator. The
second term has degree $n+1$ in $\om$ and gives rise to a falling down of the characteristic function, since the coefficient has negative sign.
So for high frequencies the filter will be a low pass filter
and for low frequencies the filter works well for differentiation.
Increase of $n$ brings the filter in a sense
closer to the ideal differentiator (see also Remark \ref{147}
for approximation to $f^{(m)}(x)$ in the $x$-domain),
but the pass band will also increase, causing
more high frequency noise
(see section 5.2 and Figure \ref{fig:1} for an example in the Legendre case).
In the practice of the construction of a differentiating filter one has to
decide to what frequency the differentiation must do the job
and how much noise one accepts. This all depends on the frequency contents
of the signal and the noise.
See also the discussion for the case of
constant weights by Barak \cite[p.2761]{39} (for $m=0$)
and by Luo et al.\ \cite[\S5]{40} (for general $m$).
\end{remark}
\subsection{Smoothing filters}
For $m=0$ the filters given by the \LHS s of \eqref{29} and \eqref{79} are
examples of {\em smoothing filters}. These have a very long
history, see Schoenberg \cite{42}, \cite{41} and references given there,
which go back as far as De Forest's work in 1878.
We put $M:=-N$ in \eqref{123} and $M=-N=-1$ in \eqref{124}:
\begin{align}
g(y)&=\sum_{x=-N}^N f(y-x)\,\rho(x)\qquad(y\in\ZZ),
\label{83}\\
g(y)&=\int_{-1}^1 f(y-x)\,\rho(x)\,dx\qquad(y\in\RR),
\label{84}
\end{align}
and we usually take $\rho$ symmetric: $\rho(x)=\rho(-x)$.

We say that \eqref{83} or \eqref{84} is {\em exact} for the degree $j$
(where $j<2N$ in case of \eqref{83}) if $g=f$ whenever
$f$ is a polynomial of degree $\le j$, but $g\ne f$ for some polynomial
$f$ of degree $j+1$. Because of symmetry of $\rho$,
such $j$ will always be odd.

Exactness of \eqref{83} for degree at least $2n+1<2N$ can equivalently
be stated as
\begin{equation*}
\rho(x)=\sum_{k=0}^{N} c_k\,t_{2k}(x+N,2N+1)\qquad(x\in\{-N,-N+1,\ldots,N\})
\end{equation*}
with
\begin{equation*}
c_k=t_{2k}(N,2N+1)/h_{2k}\quad\mbox{if $k=0,1,\ldots,n$,}
\end{equation*}
where $y\mapsto t_{2k}(y,2N+1)$ is a Gram polynomial (see \eqref{11}) and
$h_{2k}$ is the corresponding constant given by \eqref{10}.
For such $\rho$ the sum of squares
$\sum_{x=-N}^N \rho(x)^2$
is minimal if and only if $c_k=0$ for $n<k\le N$. Then
\begin{equation}
\rho(x)={\bf K}_{2n}(x+N,N)\qquad(x\in\{-N,-N+1,\ldots,N\}),
\label{89}
\end{equation}
where ${\bf K}_{2n}$ is the Christoffel-Darboux kernel of degree $2n$
(see \eqref{31})
for the orthogonal polynomials $y\mapsto t_{k}(y,2N+1)$.

Similarly, in case of \eqref{84} and assuming that $\rho$ is a polynomial,
the requirements that the formula is exact for the degree $2n+1$ and that
$\int_{-1}^1 \rho(x)^2\,dx$ is minimal are equivalent to
\begin{equation}
\rho(x)={\bf K}_{2n}(x,0)\qquad(x\in[-1,1]),
\label{91}
\end{equation}
where ${\bf K}_{2n}$ is
the Christoffel-Darboux kernel of degree $2n$ for the Legendre polynomials
$P_k$.

More generally than \eqref{89} we can work with orthogonal polynomials
$p_n$ satisfying \eqref{87} or \eqref{88} for equidistant points running
over a symmetric set and with symmetric weights:
\begin{equation*}
\sum_{x=-N}^N p_m(x)\,p_n(x) w(x)=0\qquad(m\ne n),
\end{equation*}
where $w(x)=w(-x)$.
In terms of these polynomials $p_n$
exactness of \eqref{83} for degree at least $2n+1<2N$ can equivalently
be stated as
\begin{equation}
\rho(x)=\sum_{k=0}^{N} c_k\,p_{2k}(x)\,w(x)\qquad(x\in\{-N,-N+1,\ldots,N\})
\label{92}
\end{equation}
with
\begin{equation}
c_k=p_{2k}(0)/h_{2k}\quad\mbox{if $k=0,1,\ldots,n$.}
\label{93}
\end{equation}
In particular, the choice
\begin{equation}
\rho(x):={\bf K}_{2n}(x,0)\,w(x)=
\frac{k_{2n}\,p_{2n}(0)}{k_{2n+1}h_{2n}}\,\frac{p_{2n+1}(x)\,w(x)}x\,,
\label{90}
\end{equation}
with ${\bf K}_{2n}$ the Christoffel-Darboux kernel and with \eqref{62}
used for the second equality, will make \eqref{83}
exact for the degree $2n+1$.

The {\em characteristic function} for \eqref{83} respectively \eqref{84}
is defined by \eqref{121} with $M=-N$ respectively \eqref{122} with
$M=-N=-1$.
The condition that \eqref{83} or \eqref{84} is exact for the degree $2n+1$
is equivalent with the condition that $\phi$ has power series of the form
\begin{equation}
\phi(\om)=1-a\om^{2n+2}+\ldots
\label{85}
\end{equation}
with $a\ne0$.

An $m$-fold iteration of \eqref{83} (with $N=\iy$ for convenience) yields
\begin{equation*}
g_m(y)=\sum_{x=-\iy}^\iy f(y-x)\,\rho_m(x)\qquad(y\in\ZZ),
\end{equation*}
where $\rho_m=\rho*\ldots*\rho$ is an $m$-fold convolution product.
De Forest (1878) raised the question for which choices of $\rho$ the
asymptotic behaviour of $\rho_m(x)$ for large $m$ can be described.
Schoenberg \cite[Theorem 1 and Remark 1 on p.358]{42} showed that this
is possible precisely if the characteristic function satisfies
\begin{equation}
|\phi(\om)|<1 \quad\mbox{for $0<\om<2\pi$.}
\label{86}
\end{equation}
If \eqref{86} is satisfied then the smoothing is called {\em stable}.
Clearly \eqref{86} will imply that \eqref{85} holds with $a>0$.
For $\rho$ given by \eqref{90} we see from \eqref{80} that $a>0$ is
satisfied for any choice of the weights and
that $a$ is explicitly given by
\begin{equation*}
a=\frac{|p_{2n+2}(0)|}{|k_{2n+2}|\,(2n+2)!}\,.
\end{equation*}

The stability condition \eqref{86} can also be considered for the
continuous smoothing formula \eqref{84}, where now $\om\in\RR\backslash\{0\}$
in \eqref{86}. For the Legendre case where $\rho$ is given by \eqref{91},
stability was proved by Trench \cite{43} and Lorch \& Szego \cite{44}.
\subsection{Fourier-Bessel functions}
As a common generalization of \eqref{83} and \eqref{84} with $\rho$
given by \eqref{91} and \eqref{90}, respectively, we can
consider a smoothing formula
\begin{equation}
g(y)=\int_\RR f(y-x)\,r(x)\,d\mu(x)
\end{equation}
with $\mu$ an even positive orthogonality measure
for the orthogonal polynomials
$p_n$ and with $r$ given by
\begin{equation}
r(x):={\bf K}_{2n}(x,0)=\sum_{j=0}^n\frac{p_{2j}(0)\,p_{2j}(x)}{h_{2j}}
=\frac{k_{2n}\,p_{2n}(0)}{k_{2n+1}h_{2n}}\,\frac{p_{2n+1}(x)}x\,.
\label{140}
\end{equation}
(Such usage of the Christoffel-Darboux formula was emphasized in \cite{27}
for the cases of Gram and Legendre polynomials.)$\;$
Then we can define the corresponding characteristic function by
\begin{align}
\phi(\om)&:=\int_\RR r(x)\,e^{-ix\om}\,d\mu(x)
\label{141}\\
&\;=\sum_{j=0}^n\frac{p_{2j}(0)}{h_{2j}}\,
\int_\RR p_{2j}(x)\,e^{-ix\om}\,d\mu(x)\nonumber\\
&\;=\frac{k_{2n}\,p_{2n}(0)}{k_{2n+1}h_{2n}}\,
\int_\RR \frac{p_{2n+1}(x)}x\,e^{-ix\om}\,d\mu(x).\nonumber
\end{align}
Thus
\begin{equation}
\phi'(\om)=\frac{\,k_{2n}\,p_{2n}(0)}{i\,k_{2n+1}h_{2n}}\,
\int_\RR p_{2n+1}(x)\,e^{-ix\om}\,d\mu(x).
\label{139}
\end{equation}
In the Legendre case $d\mu(x)=dx$ with support $[-1,1]$ we can
evaluate integrals occurring above in terms of (spherical)
Bessel functions as follows:
\begin{equation}
\int_{-1}^1 P_n(x)\,e^{-ix\om}\,dx=i^{-n}\sqrt{\frac{2\pi}\om}J_{n+\half}(\om)
=2i^{-n}\,j_n(\om),
\label{102}
\end{equation}
see \cite[(18.17.19), (10.47.3)]{53} or \cite[(4)]{23}.
In the Chebyshev case $d\mu(x)=(1-x^2)^{-1/2}\,dx$,
$T_n(\cos\tha):=\cos(n\tha)$  we similarly obtain
(see \cite[(10.9.2)]{53})
\begin{equation}
\pi^{-1}\int_{-1}^1 T_n(x)\,e^{-ix\om}\,(1-x^2)^{-1/2}\,dx=i^{-n}\,J_n(\om).
\label{138}
\end{equation}
Formula \eqref{138} was the reason for Mantica \cite{54},
\cite{62} to call the functions
\begin{equation}
\FSJ_n(\om;\mu):=\int_\RR p_n(x)\,e^{-ix\om}\,d\mu(x)
\label{128}
\end{equation}
{\em Fourier-Bessel functions}
(where he took $p_n$ orthonormal and $\mu$ a probability measure).
The same functions occur in Ignjatovi\'c \cite{61} as the functions
$\FSK^n[m]$ (in the notation of \cite[\S2.1]{61}).
Formula \eqref{102} played an important role in the proof of the
stability result in the Legendre case, see \cite{44}.
It also occurred in Rangarajan et al.\ \cite[(14), (15)]{24}
for a formal operational
calculus in connection with the \RHS\ of \eqref{73} before taking limits.
In the Appendix we will compute the Fourier-Bessel functions
for the case of the shifted symmetric Hahn polynomials \eqref{129}.
\subsection{Stability of smoothing in case of symmetric Hahn and
Krawtchouk polynomials}
The shifted symmetric Hahn polynomials
\begin{equation}
p_n(x):=Q_n(N+x;\al,\al,2N)\qquad(n=0,1,\ldots,2N)
\label{129}
\end{equation}
are orthogonal polynomials on $\{-N,N+1,\ldots,N\}$
with respect to the symmetric weights
\begin{equation}
w_x:=
\frac{(\al+1)_{N+x}}{(N+x)!}\,\frac{(\al+1)_{N-x}}{(N-x)!}\,,
\label{95}
\end{equation}
see \cite[(9.5.1) and (9.5.2)]{5}. Assume that $\al$ is a nonnegative
integer. Consider in terms of these polynomials $p_n$ formulas
\eqref{92} and \eqref{93} characterizing exactness for degree at least $2n+1$.
Now observe that, by \cite[(9.5.9)]{5}, we have
\begin{equation*}
\De_x^\al(w_x\,p_n(x))=\frac{(2N+1)_\al}{\al!}\,
Q_{n+\al}(N+x+\al;0,0,2N+\al),
\end{equation*}
where $\De_x(f(x))=(\De f)(x):=f(x+1)-f(x)$. It follows that
\begin{equation}
\sum_{x=-N-\al}^N ((\De^\al\rho)(x))^2
\label{94}
\end{equation}
is minimal for $\rho$ given by \eqref{92} and \eqref{93} if and only if
$c_k=0$ for $n<k\le N$.
Greville \cite[\S3]{45} (1966)
denotes \eqref{94} by $R_\al^2$ (after division
by $\binom{2\al}\al$). Therefore he calls the smoothing formula \eqref{83}
the {\em minimum $R_\al$ formula} if $\rho$ is taken such that \eqref{94}
is minimal. Greville \cite[(4.2)]{45}
gives an explicit formula for the characteristic
function $\phi$ in case of a minimum $R_\al$ formula.
He ascribes this formula to Sheppard \cite{63} (1913).
We will derive this formula in the Appendix.
Greville  \cite[\S5]{45} next proves the
stability property \eqref{86} for these cases.
Curiously, Greville does not
mention Hahn polynomials in any way. Hahn polynomials in this
context seem to come up first in Bromba \& Ziegler \cite[\S3.2]{46}.

As a special case of \cite[(9.5.16)]{5} there is the
limit formula
\begin{equation}
\lim_{\al\to\iy}Q_n(x+N;\al,\al,2N)=K_n(x+N;\thalf,2N)
=\hyp21{-n,-N-x}{-2N}2
\label{96}
\end{equation}
where the polynomials $x\mapsto K_n(x;p,N)$ are
{\em Krawtchouk polynomials} (see \cite[\S9.11]{5}).
The corresponding weights \eqref{95}, suitably normalized,
tend for $\al\to\iy$ to the symmetric weights
\begin{equation*}
w_x:=\binom{2N}{N+x}\qquad(x=-N,-N+1,\ldots,N),
\end{equation*}
with respect to which the polynomials $p_n(x)=K_n(x+N;\thalf,2N)$
are orthogonal. The smoothing formula with $\rho$ given by \eqref{90}
for this $p_n$ and $w$ is called the {\em minimum $R_\iy$ formula}
by Greville  \cite[\S6]{45}. He obtains the characteristic function $\phi$
for this case explicitly as a limit case of his formula in the
minimum $R_\al$ case (see also \eqref{146}),
not working with Krawtchouk polynomials at all.
(Krawtchouk polynomials seem to come up first in this
context in Bromba \& Ziegler \cite[\S3.3]{46}.)$\;$
But Greville also obtains in some way that,
\begin{equation}
\phi^{(k)}(\pi)=0\qquad(k=0,1,\ldots,2N-2n-1).
\label{97}
\end{equation}
Let us prove this by observing from \eqref{96} that
\[
p_{2N}(x)=\hyp21{-2N,-N-x}{-2N}2=
\sum_{j=0}^{N+x}\binom{N+x}j(-2)^j=(-1)^{N+x}.
\]
Hence
\[
\phi^{(k)}(\pi)=\sum_{x=-N}^N \rho(x)\,(-ix)^k\,(-1)^x
=i^{2N-k}\sum_{x=-N}^N {\bf K}_{2n}(x,0)\,x^k\,p_{2N}(x)\,w_x=0
\]
for $k<2N-2n$ by orthogonality. Now Greville concludes from
\eqref{85} and \eqref{97} that
\begin{equation}
\phi(\om)=1-(\sin^2(\om/2))^{n+1}P(\sin^2(\om/2))=
(\cos^2(\om/2))^{N-n}Q(\sin^2(\om/2))
\label{98}
\end{equation}
for certain polynomials $P$ of degree $N-n-1$ and $Q$ of degree $n$.
From that he immediately derives that $Q(z)$ is equal to
the power series of $(1-z)^{-N+n}$ truncated after the term with $z^n$, i.e.,
\begin{equation}
Q(z)=\sum_{k=0}^n\frac{(N-n)_k}{k!}\,z^k>0\qquad(0\le z\le 1).
\label{99}
\end{equation}
By a similar argument we see that
\begin{equation}
P(z)=\sum_{k=0}^{N-n-1}\frac{(n+1)_k}{k!}\,(1-z)^k>0\qquad(0\le z\le 1).
\label{100}
\end{equation}
Hence, by \eqref{98},
\begin{equation*}
0\le\phi(\om)<1\qquad(0<\om<2\pi),
\end{equation*}
which is even stronger than the stability condition \eqref{86}.

Identity \eqref{98} with $Q$ and $P$ given by \eqref{99}, \eqref{100}
has a long history which is surveyed in
Koornwinder \& Schlosser \cite{47}. However, this paper missed Greville's
paper and the connection with Krawtchouk polynomials.
A sequel \cite{71} to \cite{47}, tracing \eqref{98}--\eqref{100}
back to 1713, has appeared.

By a short chain of identities, using \eqref{98} and \eqref{99},
$\phi(\om)$ can be expressed as
\begin{equation*}
\phi(\om)=(N-n)\,\binom Nn\,\int_0^{\cos^2(\om/2)} s^{N-n-1}(1-s)^n\,ds,
\end{equation*}
see \cite[top of p.249]{47}. Hence $\phi$ is monotonically decreasing
on $[0,\pi]$ from 1 to 0. Such filters without ripples are called
{\em maximally flat} by Herrmann \cite{48}
(see also Samadi and Nishihara \cite{66}). Herrmann gave the same
argument as above for solving \eqref{98}, apparently unaware
of Greville~\cite{45}.
\subsection{The Savitzky-Golay paper and its follow-up}
The first instance of an approximation of first and higher derivatives
by formula \eqref{29} with possibly $n>m$ was given by
Savitzky \& Golay \cite{10} in 1964. They only dealt with the case
of constant weights on $\{-N,-N+1,\ldots,N\}$, they used only very special
$N$, $n$ and $m$, and they did not explicitly mention or use the
corresponding orthogonal polynomials. They were motivated by applications
in spectroscopy. Their paper had an enormous impact, for instance
5432 citations in Google Scholar in January 2012.
Some corrections to \cite{10} were given by
Steinier et al.\ \cite{49} in 1972.

Probably, Gorry \cite{50} (1990) was the first who gave \eqref{29} in a more
structural form in the case of constant weights on an equidistant set
using centered Gram polynomials. Next, in \cite{51} (1991) he considered
\eqref{29} on a finite non-equidistant set, still with constant weights.
Meer \& Weiss \cite{8} (1992) gave \eqref{29}
for orthogonal polynomials on a set
$\{-N,-N+1,\ldots,N\}$ with respect to general weights.
They made this more explicit in the cases
of centered Gram polynomials and of centered Krawtchouk polynomials
with symmetric weights.

We recommend Luo et al.\ \cite{40} as a relatively recent survey of
follow-up to the Savitzky-Golay paper.
\section{Filter properties in the frequency domain: some examples}
In section \ref{137} many so-called {\em linear} filters for derivatives
were mentioned.
In electrical engineering one uses {\em the transfer function} $H$
for understanding the properties of the filter.
This is the Fourier transform of the unit impulse response of
the filter, see \eqref{121}, \eqref{122} where we wrote $\phi$ instead of
$H$.
In general, this function is complex-valued. We can show
the properties of the filter in the frequency domain
by a log-log plot of the modulus of the
transfer function (which may be complemented by a phase plot).
In general, for a differentiation filter of order $n$, $H(\om)$
should behave for low frequencies like $\om^n$,
and for high frequencies like a constant (equal to zero in the ideal case).
When the behaviour is different, the filter is called {\em unstable}.

\subsection{The Lanczos derivative}
The (analog) filter corresponding to the Lanczos derivative is given by
\eqref{69} ignoring the limit:
\begin{equation}
g(x)=\frac{3}{2\delta}\,\int_{-1}^{1}f(x+\xi\delta)\,\xi\,d\xi.
\label{112}
\end{equation}
The output function $g$ can be considered as a continuous (i.e.\ unsampled)
approximation of the first derivative of the input function $f$.
The transfer function $H(\om)$ is equal to the quotient of $g$ and $f$
with $f(y):=e^{i\om y}$. A short computation gives
\begin{align}
H(\om)&=\frac3{2\de}\,\int_{-1}^1 e^{i\de\om\xi}\,\xi\,d\xi
=i\om\,\frac3{(\de\om)^3}\,\Big(\sin(\de\om)-\de\om\,\cos(\de\om)\Big)
\label{132}\\
&=i\om(1+\FSO(\de^2\om^2))\quad\mbox{as $\de\om\downarrow0$,}
\nonumber
\end{align}
compatible with \eqref{80} for $m=1$, $n=2$. For small $\om\de$ we have
$\arg(H(\om))=\pi/2$.
The modulus of $H(\om)$ for $\delta=1$ is given in Figure \ref{fig:1},
case $n=1$ as a log-log plot.

\subsection{Multi-term variant of the Lanczos derivative}
\begin{figure}[ht]
\centering     
\includegraphics[width=7cm]{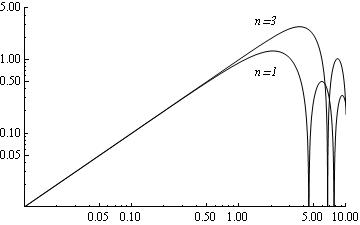} 
\caption{Modulus of transfer function for the first order analog filter,
$n=1$ and 3}
\label{fig:1}
\end{figure}

To get a better approximation we can use equation \eqref{29}
with $p_n$ a Legendre polynomial:
\begin{equation}
g(x)=\frac1{\delta^m}\,\sum_{k=0}^{[(n-m)/2]}
\frac{P_{m+2k}^{(m)}(0)}{h_{m+2k}}
\int_{-1}^1 f(x+\de\xi)\,P_{m+2k}(\xi)\,d\xi.
\label{127}
\end{equation}
Here
\begin{equation*}
\frac{P_{m+2k}^{(m)}(0)}{h_{m+2k}}=\frac{2^m(-1)^k(\thalf)_{m+k}
(m+2k+\thalf)}{k!}
\end{equation*}
by \cite[10.10(26), 10.9(19)]{4} and \eqref{13}.
For $m=n=1$ \eqref{127} reduces to \eqref{112}.
For the transfer function we obtain (see also \eqref{80}):
\begin{align*}
H_{m,n}(\om)&=\frac{2^m}{\delta^m}\,\sum_{k=0}^{[(n-m)/2]}
\frac{(-1)^k(\thalf)_{m+k}(m+2k+\thalf)}{k!}\,
\int_{-1}^{1}\,e^{i \om \delta \xi}P_{m+2k}(\xi)\,d\xi\\
&=\frac{2^{m+1}i^m}{\de^m}\,\sum_{k=0}^{[(n-m)/2]}
\frac{(\thalf)_{m+k}(m+2k+\thalf)}{k!}\,j_{m+2k}(\de\om),
\end{align*}
where the spherical Bessel functions $j_{m+2k}$ entered by \eqref{102}.
An explicit formula for spherical Bessel functions is given in
\cite[(10.49.2)]{53}.
In particular,
\begin{align*}
j_1(z)&=\frac{1}{z^2}\Big(\sin{z}-z\cos{z}\Big),\\
j_3(z)&=\frac{1}{z^4}\Big((15-6z^2)\sin{z}-(15-z^2)z\cos{z}\Big).
\end{align*}
Thus $H_{1,1}(\om)$ is given by \eqref{132}. After some computation
we obtain
\begin{equation}
H_{1,3}(\omega)=\frac{15i\om}2\,
\frac{(21-8(\de\om)^2)\sin(\de\om)+(-21\de\om+(\de\om)^3)\cos(\de\om)}
{(\de\om)^5}\,.
\label{125}
\end{equation}
The modulus of $H_{1,3}(\om)$ for $\delta=1$ is also given in figure
\ref{fig:1}, case
$n=3$. It is clear that for $n=3$ the plot stays close to a straight line
until higher values of $\om$ than for $n=1$.

\subsection{First order Savitzky-Golay filter}
When the input signal of the filter is given by a vector of (for convenience)
odd dimension $2N+1$ obtained by sampling a function $f$ on equidistant
points $x-N\de, x-(N-1)\de,\ldots,x+N\de$, then we may use \eqref{74}
ignoring the limit as a discrete
filter for the $n$-th derivative of $f$ at $x$.
In particular, for $n=1$, we can use \eqref{70}:
\[
g(x)=\frac3{2N(N+\thalf)(N+1)\de}\,
\sum_{\xi=-N}^N f(x+\de\xi)\,\xi.
\]
For the transfer function
\begin{equation}
H(\om)=\frac3{2N(N+\thalf)(N+1)\de}\,
\sum_{\xi=-N}^N e^{i\de\om\xi}\,\xi
\label{150}
\end{equation}
we obtain by straightforward computation that
\begin{equation}
H(\om)=\frac{3i}{2(2N+1)\de}\,\big(\sin(\thalf\de\om)\big)^{-2}\,
\Big(\frac{\sin(N\de\om)}N-\frac{\sin((N+1)\de\om)}{N+1}\,\Big).
\label{151}
\end{equation}
Note that the phase shift is exactly $\pi/2$.

\begin{figure}[ht]
   \centering     
\includegraphics[width=7cm]{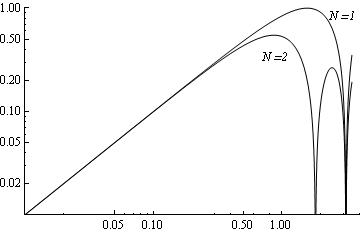} 
\caption{Modulus of transfer function for the first order discrete filter,
$N=1$ and 2}
\label{fig:2}
\end{figure}

The modulus of $H(\om)$ for $\de=1$ and for $N=1$ and 2 is given in
Figure \ref{fig:2}. See \cite[Figure 1]{40} for similar pictures.

\subsection{Butterworth filter}
If one needs a filter that does differentiation for low frequencies
very well and has a good suppression for high frequencies then there
are better filters than the ones discussed in this paper until
here. For example there are the so-called Tchebyshev, inverse
Tchebyshev, Elliptic, Butterworth and Bessel filters. These filters
all differentiate, but the choice of the most
suitable filter depends on the properties one needs, for instance a
constant phase response, a good amplitude response, less side-lobes
etc. We mention here the so-called $n$-th order
{\em Butterworth filter}.
The square of the modulus of the transfer function of an
$n$-th order analog Butterworth filter that differentiates with order $m$ is
given by
\begin{equation}
\Big|{H_{m,n}(\om)\Big|}^2=\frac{\om^{2m}}{1+(\om/\om_0)^{2n}}
=\om^{2m}\,|H_{0,n}(\om)|^2\qquad
(n>m).
\label{149}
\end{equation}
Here $\om_0$ is the so-called {\em cutoff frequency}.
It is at this frequency $\om_0$ where the
the asymptotics of the low frequency part and the high frequency part of
the transfer function meet.
The factor $\om^{2m}$ is the square of the modulus of the transfer
function $(i\om)^m$ of the ideal $m$-th order differentiator.

\begin{figure}[ht]
   \centering     
\includegraphics[width=7cm]{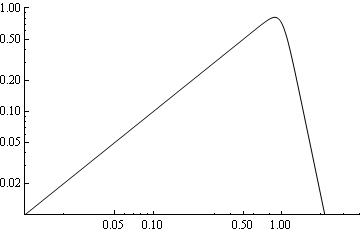} 
\caption{Transfer function for seventh order Butterworth filter}
\label{fig:3}
\end{figure}

As an example see Figure \ref{fig:3} showing the transfer function of
a seventh order Butterworth filter with $\om_0=1$
(see how the side lobes differ from those of
Figure \ref{fig:2}).

In \eqref{149} one has to make a choice of $H_{m,n}(\om)$ as follows:
\begin{equation*}
H_{m,n}(\om)=\frac{(i\om)^m}{p_n(i\om)}
\quad\mbox{with $|p_n(i\om)|^2=1+(\om/\om_0)^{2n}$}
\end{equation*}
such that $p_n$ is a polynomial of degree $n$ with real coefficients
for which all (possibly complex) roots have negative real part.
Then \eqref{122} and \eqref{124} take the form
\begin{equation*}
H_{0,n}(\om)=\int_0^\iy \rho(t)\,e^{-i\om t}\,dt,\qquad
g(t)=\int_0^\iy f(t-\tau)\,\rho(\tau)\,d\tau.
\end{equation*}
$\rho(t)$ is called the {\em impulse response} of the filter.
It follows that the output function $g$ satisfies a differential equation
with the input function $f$ as inhomogeneous part:
\[
p_n(d/dt)\,g(t)=f(t).
\]
For instance, for $n=1$ we have
\[
p_1(i\om)=1+i\om/\om_0,\quad
\rho(t)=\om_0\,e^{-\om_0t},\quad
\om_0^{-1}\,g'(t)+g(t)=f(t).
\]
For $n=2$ we have
\begin{equation*}
p_2(i\om)=1+2^{1/2}\,(i\om/\om_0)+(i\om/\om_0)^2,\quad
\rho(t)=2^{1/2}\,\om_0\,e^{-2^{-1/2}\om_0t}\,\sin(2^{-1/2}\om_0t).
\end{equation*}

One can obtain the transfer function for the Butterworth filter in the
digital case from $H_{0,n}(\om)$ in the analog case by so-called
frequency warping: replace $\om$ by $2T^{-1}\tan(\om T/2)$, where $T$
is the length of the sampling interval.
Then some linear combination of finitely many output values
$g(x),g(x-T),\ldots$ will be equal to some linear combination of
finitely many input values $f(x),f(x-T),\ldots$ (a so-called recursive
filter).

There are important differences for practical applications between
filters obtained from orthogonal polynomials, as amply considered in this
paper, and the Butterworth filter.
In the analog case the Butterworth filter can be much easier
constructed physically. But in the discrete case the filters obtained from
orthogonal polynomials are much easier to handle in the time domain
than the Butterworth filter.

For more details about the Butterworth filter for $m=0$ see
Oppenheim \& Schafer \cite[\S5.2.1]{67} (both analog and digital),
Johnson \cite[\S3.2]{29} (analog) and Hamming \cite[\S12.6]{28} (digital).
See Herely \& Vetterli \cite[\S IV-A]{68} and Cohen \& Daubechies
\cite[\S6.2]{69} for Butterworth filters in connection with wavelets.
\appendix
\section{Appendix}
Here we derive an explicit expression for the Fourier-Bessel
functions (see \eqref{128})
associated with the shifted symmetric Hahn polynomials (see \eqref{129}).
First observe that
by the duality between Hahn polynomials and dual Hahn polynomials (see the
formula after (9.6.16) in \cite{5}) the generating function
\cite[(9.6.12)]{5} for dual Hahn polynomials can be rewritten in
terms of Hahn polynomials:
\begin{equation}
(1-t)^n\,\hyp21{n-N,n+\al+1}{-\be-N}t=\frac{N!}{(\be+1)_N}\,
\sum_{x=0}^N w_x\,Q_n(x;\al,\be,N)\,t^x,
\label{130}
\end{equation}
where
\begin{equation}
w_x:=\binom{\al+x}x \binom{\be+N-x}{N-x},
\label{131}
\end{equation}
i.e., the weight occurring in the orthogonality relation
\cite[(9.5.2)]{5} for Hahn polynomials.

Next, in \eqref{130} take $\be=\al$, $t:=e^{-i\tha}$, replace $N$ by $2N$,
shift $x$ to $x+N$ and apply Pfaff's identity \cite[(2.3.14)]{34}. Then
\begin{equation}
\frac{(2\al+2)_{2N+n}}{2^{2n}(2N)!(\al+\tfrac32)_n}\,e^{iN\tha}
(1-e^{-i\tha})^n\,
\hyp21{n-2N,n+\al+1}{2n+2\al+2}{1-e^{-i\tha}}
=\sum_{x=-N}^N w_x\,p_n(x)\,e^{-ix\tha},
\label{133}
\end{equation}
where $p_n(x)$ and $w_x$ are given by \eqref{129} and \eqref{95},
respectively.

Now use the quadratic transformation \cite[2.11(30)]{60} and the expression
\cite[10.9(20)]{4} of Gegenbauer polynomials in terms of
hypergeometric functions:
\begin{align*}
\hyp21{n-2N,n+\al+1}{2n+2\al+2}{1-e^{-i\tha}}
&=e^{i(\half n-N)\tha}\,\hyp21{n-2N,2N+n+2\al+2}{n+\al+\tfrac32}
{\thalf\big(1-\cos(\thalf\tha)\big)}\\
&=\frac{(2N-n)!}{(2n+2\al+2)_{2N-n}}\,e^{i(\half n-N)\tha}\,
C_{2N-n}^{\,n+\al+1}\big(\cos(\thalf\tha)\big).
\end{align*}
Thus we can rewrite \eqref{133} as
\begin{equation}
\frac{(\al+1)_n}{(-2N)_n}\,i^{-n}\,\big(2\sin(\thalf\tha)\big)^n\,
C_{2N-n}^{\,n+\al+1}\big(\cos(\thalf\tha)\big)
=\sum_{x=-N}^N w_x\,p_n(x)\,e^{-ix\tha}.
\label{134}
\end{equation}
The \LHS\ of \eqref{134} gives an expression for the
Fourier-Bessel function \eqref{128} associated with the shifted symmetric Hahn
polynomials \eqref{129}.

Now let $\phi$ be defined by \eqref{141}, \eqref{140} with $p_n$ given
by \eqref{129}. Then
combination of \eqref{139}, \eqref{134} and \cite[10.9(22)]{4} yields that
\begin{equation}
\frac{d\phi(\om)}{d(\sin^2(\thalf\om))}=C\,(\sin^2(\thalf\om))^n\,
\hyp21{-N+n+1,N+n+\al+2}{2n+\al+\tfrac52}{\sin^2(\thalf\om)}.
\label{142}
\end{equation}
where
\begin{align}
C&=(-1)^{n+1}\,2^{2n+1}\,\frac{k_{2n}\,p_{2n}(0)}{k_{2n+1}\,h_{2n}}\,
\frac{(\al+1)_{2n+1}}{(-2N)_{2n+1}}\,
\frac{(4n+2\al+4)_{2N-2n-1}}{(2N-2n-1)!}
\nonumber\\
&=(-1)^n\frac{(N+\al+1)_{n+1}\,(-N)_{n+1}}{(n+\al+\tfrac32)_{n+1}\,n!}\,.
\label{143}
\end{align}
In the last equality we used \cite[(9.5.2), (9.5.4)]{5} for $h_n$ and $k_n$
and \cite[(2.4)]{65} together with \cite[(9.5.3)]{5} for getting
\begin{equation}
p_{2n}(0)=\frac{(\thalf)_n\,(N+\al+1)_n}{(-N+\thalf)_n\,(\al+1)_n}\,.
\label{144}
\end{equation}
Formulas \eqref{142}, \eqref{143} coincide with formulas
(4.3), (4.4) in Greville \cite{45} if we replace
$N,n,\al$ by $n,k,m$, respectively.
Integration of \eqref{143} together with $\phi(0)=1$ yields
\begin{equation}
\phi(\om)=1+\frac{(-1)^n}{n!}\,
\sum_{k=n+1}^N\frac{(N+\al+1)_k\,(-N)_k}{(n+\al+\tfrac32)_k\,(k-n-1)!}\,
\frac{(\sin^2(\thalf\om))^k}k\,.
\label{145}
\end{equation}
Formula \eqref{145} coincides with formula (4.2) in Greville \cite{45},
which Greville (in his earlier form (4.1)) ascribes to Sheppard \cite{63}.
However, we have only been able to
find a match of the special case $m=0$ of
\cite[(4.1)]{45} with a formula in Sheppard's paper, namely with
\cite[(64)]{63}.

In the limit for $\al\to\iy$ formula \eqref{145} becomes
\begin{equation}
\phi(\om)=1+\frac{(-1)^n}{n!}\,
\sum_{k=n+1}^N\frac{(-N)_k}{(k-n-1)!}\,
\frac{(\sin^2(\thalf\om))^k}k\,.
\label{146}
\end{equation}
This coincides with the formula after (6.1) in Greville \cite{45}, and
also with \eqref{98} combined with \eqref{100} and
\cite[(2.3.15)]{34}.

If $\al\in\Znonneg$ then the \LHS\ of \eqref{134} can be written as
a finite sum, where the number of terms is independent of $N$.
First observe that by \cite[(14.13.1), (14.3.21), (5.5.5)]{53} we have
\begin{multline}
C_n^{\la}(\cos\tha)=\frac{\Ga(2\la+1)}{2^{2\la}\Ga(\la+1)^2}\,
\frac{(2\la)_n}{(\la+1)_n}\,(\sin\tha)^{1-2\la}\,
\sum_{k=0}^\iy\frac{(1-\la)_k(n+1)_k}{(n+\la+1)_k k!}\,
\sin\big((2k+n+1)\tha\big)\\
(\la>0,\;0<\tha<\pi).
\label{135}
\end{multline}
For $\la\in\Zpos$ the above series terminates after the term with
$k=\la-1$. Hence, for $\al\in\Znonneg$ \eqref{134} takes the form
\begin{multline}
-2\,i^{n+2\al}\,\binom{2N+\al}\al\,\big(2\sin(\thalf\tha)\big)^{-n-2\al-1}\,
\sum_{k=0}^{n+\al}(-1)^k\,\binom{n+\al}k\,
\frac{(-2N-2\al-n-1)_k}{(-2N-\al)_k}\\
\times\sin\big(\thalf(2N+2\al+n-2k+1)\tha\big)
=\sum_{x=-N}^N w_x\,p_n(x)\,e^{ix\tha}.
\label{136}
\end{multline}
For $\al=0$ and $n=1$ formula \eqref{136} specializes to \eqref{151} with $H(\om)$
given by \eqref{150}. Use that $Q_1(N+x;0,0,2N)=x/N$.

\quad\\
\begin{footnotesize}
E. Diekema, Kooikersdreef 620, 7328 BS Apeldoorn, The Netherlands;
\sPP
email: {\tt e.diekema@gmail.com}
\bLP
T. H. Koornwinder, Korteweg-de Vries Institute, University of
Amsterdam,\sPP
P.O.\ Box 94248, 1090 GE Amsterdam, The Netherlands;
\sPP
email: {\tt thkmath@xs4all.nl}
\end{footnotesize}

\end{document}